\title{{\scshape{A Note on the Variance of the Square Components \\ of a Normal Multivariate 
      within a Euclidean Ball \\[1.0ex]}}}

\author{\color{Blue}
  {\sc  Filippo Palombi}\footnote{Corresponding author. E--mail: 
  {\tt filippo.palombi$@$istat.it}} \ {\it and}  %
  {\sc Simona Toti} \\[2.0ex]
  \includegraphics[angle=0.0,width=0.1\textwidth]{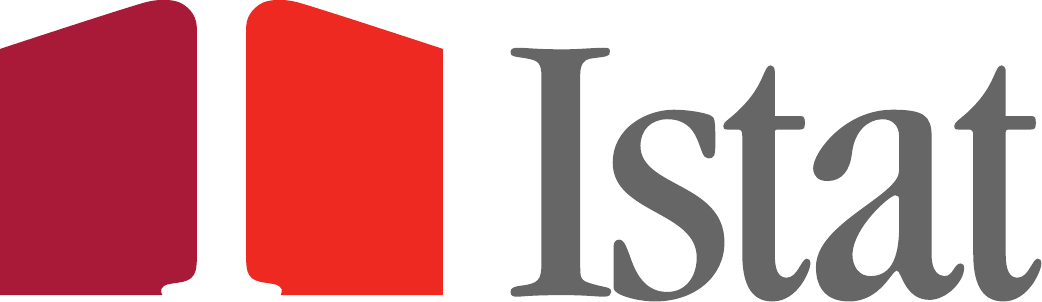}\\[2.0ex]
  \color{Blue}
  {\sc Istituto Nazionale di Statistica}\\
  \color{Blue}
  {\it Via Cesare Balbo 16, 00184 Rome -- Italy}\\[2.0ex]
}
\date{\today}
\documentclass[10pt,fleqn]{article}

\usepackage{bold-extra}
\usepackage{titlesec}
\usepackage[usenames,dvipsnames]{xcolor}
\usepackage{footmisc,booktabs}
\usepackage{psfrag}
\usepackage{amsmath}
\usepackage{amsfonts}
\usepackage{dsfont}
\usepackage{amsthm}
\usepackage{epsfig}
\usepackage{multicol}
\usepackage[title]{appendix}
\usepackage{subfigure}
\usepackage{pgf}
\usepackage{bbm}
\usepackage{tikz}
\usepackage{undertilde}
\usetikzlibrary{arrows,automata}
\usetikzlibrary{positioning}
\usepackage[utf8]{inputenc}
\usepackage[margin=30pt,font=small,labelfont=bf,labelsep=endash]{caption}
\usepackage{graphicx,xcolor}

\newenvironment{@abssec}[1]{%
     \if@twocolumn
       \section*{#1}%
     \else
       \vspace{.05in}\footnotesize
       \parindent .2in
         {\upshape\bfseries\color{Blue} #1. }\ignorespaces
     \fi}
     {\if@twocolumn\else\par\vspace{.1in}\fi}

\renewenvironment{abstract}{\begin{@abssec}{\abstractname}}{\end{@abssec}}
\renewcommand\abstractname{\color{Blue}Abstract}

\numberwithin{equation}{section}\usepackage{amssymb}
\definecolor{Black}{rgb}{0,0,0}
\definecolor{shadecolor}{rgb}{1,1.0,0.8}

\titleformat{\section}
{\color{Blue}\normalfont\Large\bfseries}
{\color{Blue}\thesection}{1em}{}

\titleformat{\subsection}
{\color{Blue}\normalfont\large\bfseries}
{\color{Blue}\thesubsection}{1em}{}

\newcommand{\lhs}{{\it l.h.s.}\ }
\newcommand{\rhs}{{\it r.h.s.}\ }

\newcommand{\uc}{\underline{c}}
\newcommand{\ud}{\underline{d}}
\newcommand{\ue}{\underline{e}}

\newcommand{\ovc}{\overline{c}}
\newcommand{\ovd}{\overline{d}}
\newcommand{\ove}{\overline{e}}

\newcommand{\cB}{{\cal B}}
\newcommand{\cC}{{\cal C}}

\newcommand{\cE}{{\cal E}}
\newcommand{\cH}{{\cal H}}

\newcommand{\cN}{{\cal N}}

\newcommand{\cP}{{\cal P}}
\newcommand{\cQ}{{\cal Q}}

\newcommand{\cS}{{\cal S}}

\newcommand{\cX}{{\cal X}}

\newcommand{\fD}{{\frak D}}

\newcommand{\fDn}{{\fD}_{\rm n}}
\newcommand{\fDd}{{\fD}_{\rm d}}

\newcommand{\rd}{{\rm d}}

\newcommand{\cov}{{\rm cov}}
\newcommand{\var}{{\rm var}}

\newcommand{\re}{{\rm e}}

\newcommand{\E}{{\mathbb{E}}}
\newcommand{\M}{{\mathbb{M}}}
\newcommand{\RR}{{\mathbb{R}}}
\newcommand{\NN}{{\mathbb{N}}}
\newcommand{\trans}[1]{{#1}^{\scriptscriptstyle{\rm T}}}

\newcommand{\sign}{{\rm sign\,}}
\newcommand{\diag}{{\rm diag}}

\newcommand{\widesim}[2][3]{
  \mathrel{\overset{#2}{\scalebox{#1}[1.0]{$\sim$}}}
}

\newtheorem{prop}{\color{Blue}Proposition}[section]
\newtheorem{conj}{\color{Blue}Conjecture}[section]

\newtheorem{theo}{\color{Blue}Theorem}[section]
\newtheorem{lemma}{\color{Blue}Lemma}[section]

\DeclareMathOperator*{\argmax}{arg\,max}

\makeatletter
\newcommand{\raisemath}[1]{\mathpalette{\raisem@th{#1}}}
\newcommand{\raisem@th}[3]{\raisebox{#1}{$#2#3$}}
\makeatother

\def\lambdabar{\protect\@lambdabar} 
\def\@lambdabar{% 
\relax 
\bgroup 
\def\@tempa{\hbox{\raise.73\ht0 
\hbox to0pt{\kern.25\wd0\vrule width.5\wd0 
height.1pt depth.1pt\hss}\box0}}% 
\mathchoice{\setbox0\hbox{$\displaystyle\lambda$}\@tempa}% 
{\setbox0\hbox{$\textstyle\lambda$}\@tempa}% 
{\setbox0\hbox{$\scriptstyle\lambda$}\@tempa}% 
{\setbox0\hbox{$\scriptscriptstyle\lambda$}\@tempa}% 
\egroup 
} 

\renewcommand{\qedsymbol}{$\boxdot$\ }

\begin{document}

\maketitle
\vskip -0.3cm
\begin{abstract}
We present arguments in favor of the inequalities
$\var(X_n^2\,|\,X\in\cB_v(\rho))\le
2\lambda_n\E[X_n^2\,|\,X\in\cB_v(\rho)]$, where
$X\sim\cN_v(0,\Lambda)$ is a normal vector in $v\ge 1$ dimensions, with 
zero mean and covariance matrix $\Lambda=\diag(\lambda)$, and $\cB_v(\rho)$
is a centered $v$--dimensional Euclidean ball of square radius
$\rho$. Such relations lie at the heart of an iterative algorithm,
proposed by Palombi et al. (2012)~\cite{palombi4} to perform a reconstruction of $\Lambda$ from
the covariance matrix of $X$ conditioned to $\cB_v(\rho)$. In the
regime of strong truncation, {\it i.e.} for $\rho\lesssim\lambda_n$,
the above inequality is easily proved, whereas it becomes harder for
$\rho\gg\lambda_n$. Here, we expand both sides in a function series
controlled by powers of $\lambda_n/\rho$ and show that the coefficient
functions of the series fulfill the inequality order by order if
$\rho$ is sufficiently large. The intermediate region remains at
present an open challenge.   
\end{abstract}

\section{Introduction}

It is intuitively clear that independent random variables develop correlations
once constrained within compact multivariate domains. Whenever the mathematical framework
rules out closed--form results, a possible approach to studying such
correlations is to focus on inequalities among expected values. 
%on the subject, and certainly a seminal paper is represented by
%ref.~\cite{joagdev}, where the concept of negative association among random
%variables has been introduced.  
As a case in point, in this paper we consider a random vector $X\sim\cN_v(0,\Lambda)$ in $v\ge 1$
dimensions, with $\Lambda = \diag(\lambda)$ and $\lambda =
\{\lambda_k\}_{k=1}^{v}$, whose probability density is cut off
sharply outside a Euclidean ball
\begin{equation} 
\cB_v(\rho) = \{x\in\RR^v\ : \ \trans{x}x <\rho\}\,.
\end{equation}
Owing to the symmetry mismatch between $\cN_v(0,\Lambda)$ and $\cB_v(\rho)$,
%for $\lambda \ne \{\tilde\lambda,\ldots,\tilde\lambda\}$ and $\tilde\lambda>0$
the conditional moments of $X$ admit no exact representation in terms of
elementary functions. Our aim is to show that the effect of the spherical
truncation on the variance of the square components of $X$ is quantified by
the 
inequalities  
\begin{equation}
\Delta_n(\rho;\lambda) \,\equiv\, \frac{1}{\rho^2}\left\{\var\left(X_n^2\,|\,X\in\cB_v(\rho)\right) -
2\lambda_n\E[X_n^2\,|\,X\in\cB_v(\rho)]\right\}\, \le 0\,,\qquad n=1,\ldots v\,.
\label{eq:varineq}
\end{equation}
The interest we have in eq.~(\ref{eq:varineq}) originates from
ref.~\cite{palombi4}, where we have proposed a fixed--point algorithm for the
reconstruction of $\Lambda$, in case the only available
information amounts to the covariance matrix ${\frak S}_\cB$ of $X$ conditioned to
$\cB_v(\rho)$. In particular, in that paper we showed that eq.~(\ref{eq:varineq}) and
the correlation inequalities 
\begin{equation}
\cov(X_n^2,X_m^2\,|\,X\in\cB_v(\rho))\le 0\,,\qquad n\ne m\,,
\label{eq:covineq}
\end{equation}
are necessary and sufficient for the convergence of the
algorithm. Eq.~(\ref{eq:covineq}) expresses a property of negative association
among the square components of $X$ (see ref.~\cite{joag:1983}). A proof of
it goes beyond the scope of the present paper. 

\begin{figure}
\centering
%\hskip 1.2cm
\includegraphics[width=0.80\textwidth]{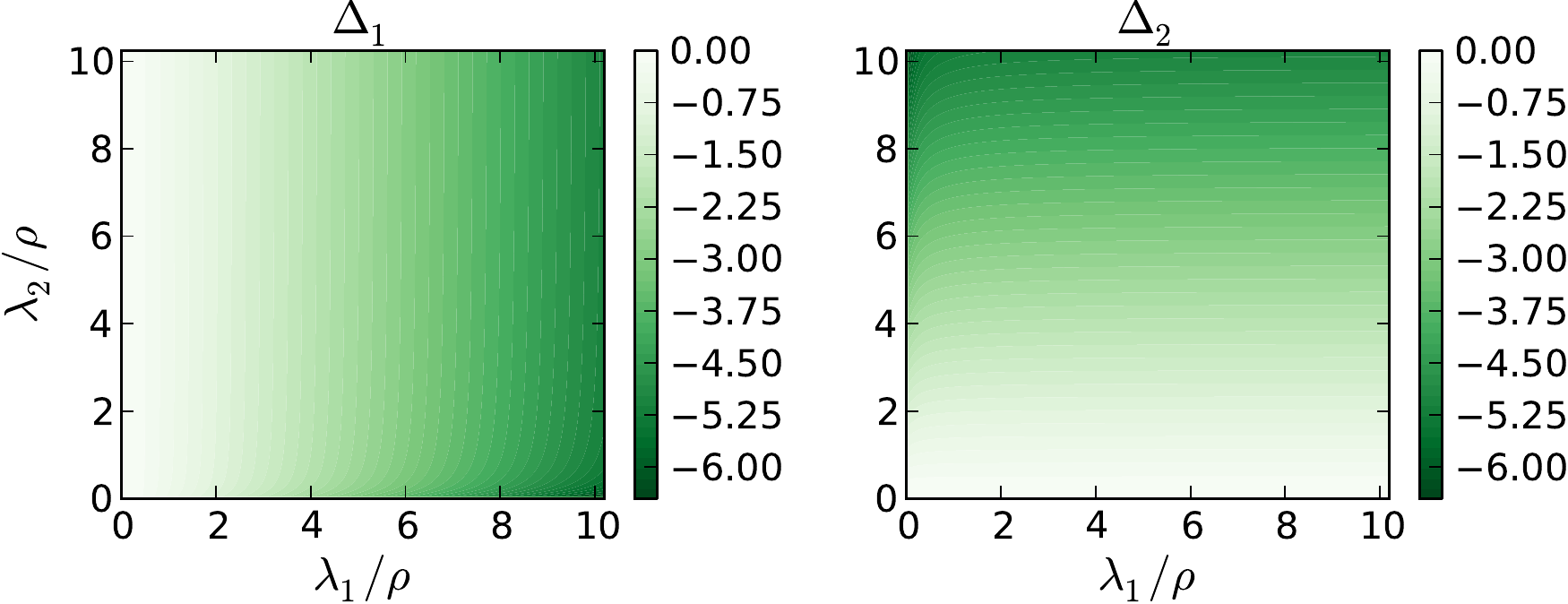}
\caption{\small Contour plots of $\Delta_n$ at $v=2$.}
\label{fig:contour}
\end{figure}

If we denote Gaussian integrals over Euclidean balls by  
\begin{equation}
\alpha_{k\ell m\dots}(\rho;\lambda) = \int_{\cB_v(\rho)}\rd^v x\
\frac{x_k^2}{\lambda_k}\,\frac{x_\ell^2}{\lambda_\ell}\,\frac{x_m^2}{\lambda_m}\dots\,\prod_{j=1}^v\delta(x_j,\lambda_j)\,,
\qquad \delta(y,\eta)\, =\,
\frac{\,\re^{-y^2/({2\eta})}\,}{(2\pi\eta)^{1/2}}\,,
\label{eq:alphaints}
\end{equation}
and define $\partial_n \equiv \partial/(\partial\lambda_n)$, we see that $\Delta_n =
2(\lambda_n^2/\rho^2)[\lambda_n\partial_n(\alpha_n/\alpha)]$. Thus, the inequality $\Delta_n \le
0$ holds true iff ${\alpha_n}/{\alpha}$ is monotonic decreasing in
$\lambda_n$ with $\lambda_{(n)} \equiv \{\lambda_i\}_{i\ne n}$ kept fixed. Since such
monotonic behavior is held by both
$\alpha_n$ and $\alpha$ separately, eq.~(\ref{eq:varineq}) simply means that
$\alpha_n$ is more rapidly decreasing than $\alpha$. An illustrative example
is shown in Fig.~\ref{fig:contour}, where contour plots of $\Delta_n$
at $v=2$ are reported\footnote{Numerical techniques for the computation of $\alpha_{k\ell
    m\ldots}$ are discussed in ref.\cite{palombi4}.}. An analysis of the monotonic properties
of averages of monotonic observables of P\'oyla distributions under linear constraints
has been originally proposed by Efron \cite{efron}. Unfortunately, the techniques there presented  
do not carry over to our set--up. 

Apart from  the covariance reconstruction algorithm, a
different motivation to care about  eqs.~(\ref{eq:varineq}) and
(\ref{eq:covineq}) has to do with non--linear optimization issues. Thanks to
Pr\'ekopa's Theorem \cite{prekopa}, $\alpha(\rho;\lambda)$ is easily shown to
be logarithmic concave in $\rho$. In sect.~2 we discuss how log--concavity
relates to the correlation inequalities. The outcome is that, independently of
Pr\'ekopa's Theorem, eqs.~(\ref{eq:varineq}) and (\ref{eq:covineq}) are alone
sufficient to induce log--concavity, while they cannot be deduced from it.

Now, despite the seemingly candid aspect of eq.~(\ref{eq:varineq}), it is 
difficult to find a unique rigorous proof of it, which works across 
the whole parameter space $(\rho;\lambda)\in\RR_+\times\RR_+^v$. In
this paper we propose three different partial arguments. The first two are discussed in 
sect.~3. They are both straightforward and apply in the  {\it
  regime of strong truncation}, {\it i.e.} for
$0<\rho<\lambda_n$ resp. $0<\rho<2\lambda_n$, independently of
$\lambda_{(n)}$. In particular, the first one is based on H\"older's
inequality, while the second one makes use of the integral representation of
$\Delta_n$. The relative ease of proving
eq.~(\ref{eq:varineq}) at strong truncation is certainly due to the large
negative values $\Delta_n$ assumes in this regime and its weak dependence upon
$\lambda_{(n)}$, which in a sense~makes~the~problem {\it   nearly}
{1--{di\-mensional}}. 

The third argument, presented in sects.~4 and 5, applies instead
in the {\it regime of weak truncation}, {\it i.e.} for $\rho\gg\lambda_n$,
where proving eq.~(\ref{eq:varineq}) is definitely harder. As
$\rho\to\infty$, we have indeed $\var(X_n^2\,|\,X\in\cB_v(\rho))\to 2\lambda_n^2$ and
$\E[X_n^2\,|\,X\in\cB_v(\rho)]\to \lambda_n$, thus $\Delta_n\to 0$. Hence, if
eq.~(\ref{eq:varineq}) is correct, it must  follow
from a cancellation of two positive terms resulting in an increasingly small
negative balance.   Motivated by the observation that the volume
constraint weakens as $\rho\to\infty$ (and consequently $\alpha_{k\ell
  m\ldots}$ becomes well approximated by a product of 1--dimensional Gaussian
integrals), we expand $\Delta_n$ in a non--elementary--function series
around the factorization point. Each term of the expansion factorizes into a
$1$--dimensional integral along the $n^{\rm th}$ direction plus a residual
$(v-1)$--dimensional integral in the orthogonal subspace. We prove that
such factors get opposite signs as $\rho\to\infty$, thus resulting in 
negative contributions.    

We finally draw our conclusions in sect.~6.

Distributional truncations find application in several frameworks. The specific one,
considered in the present paper, turns out to be useful for the compositional analysis of multivariate 
log--normal data affected by outlying contaminations, where the spherical truncation
corresponds to keeping only data with a square Aitchison distance from the mean below
a given threshold. This idea is discussed for instance in ref.~\cite{palombi1}. In that context, 
the iterative algorithm of ref.~\cite{palombi4} allows for an estimate of the complete covariance
matrix from its truncated counterpart. 

\section{Basic properties of $\alpha_{k\ell m\ldots}$}

It is worthwhile starting our discussion by reviewing some trivial 
properties of the Gaussian integrals, which are used in the sequel.
As an alternative notation for $\alpha_{k\ell m\ldots}$ we sometimes adopt the
symbol $\alpha_{1:k_1\ldots v:k_v}$, where $k_j$ counts the
multiplicity of the index $j=1,\ldots,v$. Whenever a directional index has
zero multiplicity, we simply drop it. For instance, we write
$\alpha_{j:k_j}$ in place of the more pedantic $\alpha_{1:0\ldots j:k_j\ldots
v:0}$.  
When needed, we declare the integral dimension of $\alpha_{k\ell m\ldots}$
explicitly by writing the latter as $\alpha^{(v)}_{k\ell m\ldots}$. With this in mind,
we proceed to a first set of statements.

\begin{prop}
Gaussian integrals fulfill the following properties: 
\begin{itemize}
\item[$p_1$)]{$\alpha_{k\ell m\dots}(\rho;\lambda)$ is increasing in $\rho$;}
\item[$p_2$)]{$\alpha_{k\ell m\dots}(\rho;\lambda)$ is separately decreasing
    in $\lambda_1,\ldots,\lambda_v$;}
\item[$p_3$)]{$\alpha_{k\ell m\dots}(\rho;\lambda)$ fulfills the scaling equation}
\end{itemize}
\begin{equation}
\left[\rho\partial_\rho + \sum_{r=1}^v \lambda_r\partial_r\right]
\alpha_{k\ell m\dots}(\rho;\lambda) = 0\,;
\label{eq:scaling}
\end{equation}
\begin{itemize}
\item[$p_4$)]{one--index integrals $\alpha_{k:n}$ follow the hierarchy}
\end{itemize}
\begin{equation}
\alpha_{k:n} \le (2n-1)\alpha_{k:(n-1)} \le (2n-1)(2n-3)\alpha_{k:(n-2)}\le \ldots
\le (2n-1)!!\,\alpha \le (2n-1)!!\,;
\label{eq:momrecurs}
\end{equation}
\begin{itemize}
\item[$p_5$)]{differentiating $\alpha_{k\ell m\dots}(\rho;\lambda)$ with respect
  to $\rho$ yields} 
\end{itemize}
\begin{equation}
\rho\partial_\rho \alpha_{k_1\ldots k_n}(\rho;\lambda) = \frac{1}{2}(v+2n)\alpha_{k_1\ldots
  k_n}(\rho;\lambda) -\frac{1}{2}\sum_{k=1}^v\alpha_{k_1\ldots k_n
  k}(\rho;\lambda)\,,\qquad n=0,1,2,\ldots
\label{eq:recursion}
\end{equation}
\begin{itemize}
\item[$p_6$)]{$\alpha(\rho;\lambda)$ is logarithmic concave in $\rho$,
    i.e. it fulfills}
\end{itemize}
\begin{equation}
\alpha\bigl(s\rho_1 + (1-s)\rho_2;\lambda\bigr) \ge
\bigl[\alpha(\rho_1;\lambda)\bigr]^s\bigl[\alpha(\rho_2;\lambda)\bigr]^{1-s}\,,\qquad
0\le s\le 1\,,\qquad \rho_1,\rho_2\in\RR_+\,.
\label{eq:logconcavity}
\end{equation}
\end{prop}
\begin{proof}
Property $p_1$) follows from the positiveness of the integrand of $\alpha_{k\ell
  m\ldots}$ and the observation that if $\rho_1<\rho_2$ then $\cB_v(\rho_1)\subset
\cB_v(\rho_2)$. \qedsymbol Property
$p_2$) follows from the observation that $\alpha_{k\ell m\ldots}$ depends on $\rho$ and 
$\lambda$ only via adimensional ratios, as seen by rescaling the integration variable $x\to x/\sqrt{\rho}$ in
eq.~(\ref{eq:alphaints}), {\it i.e.}
\begin{align}
\alpha_{k\ell m\dots}(\rho;\lambda) & =
\frac{\rho^{v/2}}{(2\pi)^{v/2}|\Lambda|^{1/2}} \int_{\cB_v(1)}\rd^v x\
\frac{\rho x_k^2}{\lambda_k}\frac{\rho x_\ell^2}{\lambda_\ell}\frac{\rho
  x_m^2}{\lambda_m}\cdots\, \exp\left\{-\frac{\rho}{2}\sum_{m=1}^v\frac{x_m^2}{\lambda_m}\right\} \nonumber\\[2.0ex] 
& =  \alpha_{k\ell m\ldots}\left(1;\left\{\frac{\lambda_1}{\rho},\ldots,\frac{\lambda_v}{\rho}\right\}\right) \,.
\label{eq:alpharescaled}
\end{align}
When a single variance is downscaled, e.g. $\lambda\to
\lambda'=\{\lambda_1,\ldots,a\lambda_r,\ldots,\lambda_v\}$ with $0<a<1$,  the
change in $\alpha_{k\ell m\ldots}$ is entirely transferred 
to the integration region, {\it i.e.}
\begin{align}
\alpha_{k\ell
  m\ldots}\left(\rho;\lambda'\right)
= \frac{\rho^{v/2}}{(2\pi)^{v/2}|\Lambda|^{1/2}} \int_{\cE_v(1;a)}\rd^v x\
\frac{\rho x_k^2}{\lambda_k}\frac{\rho x_\ell^2}{\lambda_\ell}\frac{\rho
  x_m^2}{\lambda_m}\cdots\, \exp\left\{-\frac{\rho}{2}\sum_{m=1}^v\frac{x_m^2}{\lambda_m}\right\}\,,
\end{align}
with
\begin{equation}
\cE_v(1;a) = \left\{x\in\mathbb{R}^v:\ x_1^2 + \dots
+ ax_r^2+\dots + x_v^2 < 1\right\}\,.
\end{equation}
Since $\cB_v(1)\subset \cE_v(1;a)$, it follows $\alpha_{k\ell m\ldots}(\rho;\lambda')> \alpha_{k\ell m\ldots}(\rho;\lambda)$. \qedsymbol
Property $p_3$) follows from the application of the chain rule of differentiation to 
eq.~(\ref{eq:alpharescaled}). The meaning of the scaling equation is that $\alpha_{k\ell
  m\ldots}$ keeps invariant under a change of the units adopted to measure both $\rho$
and $\lambda$. \qedsymbol To get convinced about property $p_4$), we first
notice that
\begin{equation}
\lambda_k\partial_k\,\alpha_{1:n_1\ldots k:n_k\ldots v:n_v} =
\frac{1}{2}\left[\alpha_{1:n_1\ldots k:(n_k+1)\ldots v:n_v} -
  (2n_k+1)\alpha_{1:n_1\ldots k:n_k\ldots v:n_v}\right]\,,
\label{eq:derlk}
\end{equation}
as proved by evaluating the derivative on the {\it l.h.s.} under the integral
sign. From property $p_2$), the {\it r.h.s.} of eq.~(\ref{eq:derlk}) is
recognized to be negative. The proof is completed by taking $n_j=0$ for $j\ne
k$. Note that eq.~(\ref{eq:momrecurs}) entails the inequalities
\begin{equation}
\E[X_k^{2n}\,|\,X\in\cB_v(\rho)]\,\le\, (2n-1)!!\, \lambda_k^n\,,\qquad n=1,2,\ldots\,,
\label{eq:moments}
\end{equation}
with the quantity on the {\it r.h.s.} representing the value of the unconditioned
 $(2n)^{\rm th}$ univariate moment of $X_k$. We shall use the
lowest order inequalities $\E[X_k^2\,|\,X\in\cB_v(\rho)]\le\lambda_k$ and
$\E[X_k^4\,|\,X\in\cB_v(\rho)]\le 3\lambda_k^2$ time and again in the
sequel. Note also that the larger $n$, the slower $\alpha_{k:n}$ saturates to its 
infinite volume limit. Indeed, if we denote by $d_n \equiv
[(2n-1)!!-\alpha_{k:n}]/\alpha_{k:n}$ the fractional distance of
$\alpha_{k:n}$ from its infinite volume limit, then eq.~(\ref{eq:momrecurs})
is equivalent to the inequality chain
\begin{equation}
{d}_0(\rho;\lambda)\le {d}_1(\rho;\lambda)\le{d}_2(\rho;\lambda)\le\ldots\,.
\label{eq:fraconv}
\end{equation}
As we shall see, this property lies at the heart of most of the difficulties
related to proving eq.~(\ref{eq:varineq}).\,\qedsymbol Property $p_5$) follows
from eqs.~(\ref{eq:scaling}) and (\ref{eq:derlk}). \qedsymbol  Finally, in
order to prove property $p_6$), we recall \cite{prekopa}

\vskip 0.2cm
\begin{theo}[Pr\'ekopa]
  \label{th:prekopa}
  Let $Q(x)$ be a convex function defined on the entire $v$--dimensional
  space $\RR^v$. Suppose that $Q(x)\ge a$, where $a$ is some real number. Let
  $\psi(z)$ be a function defined on the infinite interval $[a,\infty)$. Suppose
  that $\psi(z)$ is non--negative, non--increasing, differentiable, and
  $-\psi'(z)$ is logarithmic concave. Consider the function $f(x) = \psi(Q(x))$
  ($x\in\RR^v$) and suppose that it is a probability density, {\it i.e.}
  \begin{equation}
    \int_{\RR^v} \rd^v x\ f(x) = 1
  \end{equation}
  Denote by $P\{C\}$ the integral of $f(x)$ over the measurable subset $C$ of
  $\RR^v$. If $A$ and $B$ are any two convex sets in $\RR^v$, then the following
  inequality holds: 
  \begin{equation}
    \left(P\{A\}\right)^s
    \left(P\{B\}\right)^{1-s}\le P\{s A + (1-s)B\}\,,\qquad 0\le s \le 1\,,
  \end{equation}
  where the linear combination on the l.h.s. denotes the Minkowski sum
  \begin{equation}
    s A + (1-s)B \equiv \bigl\{s x + (1-s)y:\ x\in A,\,
    y\in B\bigr\}\,.
  \end{equation}
\end{theo}

\vskip 0.2cm
\noindent Obviously, theorem \ref{th:prekopa} applies if $f(x)$ is a product of univariate Gaussian
densities, as is the case with $\alpha(\rho;\lambda)$. In addition, if $x\in\cB_v(\rho_1)$
and $y\in\cB_v(\rho_2)$, from the convexity of the square
function $x\mapsto x^2$ it follows that
\begin{align}
\sum_{k=1}^v [sx_k + (1-s)y_k]^2 & \le s\sum_{k=1}^v x_k^2 + (1-s)\sum_{k=1}^v y_k^2 \le s\rho_1 + (1-s)\rho_2\,,
\end{align}
{\it i.e.} $s\cB_v(\rho_1) + (1-s)\cB_v(\rho_2) \subseteq \cB_v(s\rho_1 +
(1-s)\rho_2)$. Accordingly, we conclude that
\begin{align}
\left[\alpha(\rho_1;\lambda)\right]^s&\left[\alpha(\rho_2;\lambda)\right]^{1-s}
\le\ \int_{s\cB_v(\rho_1)+(1-s)\cB_v(\rho_2)}\rd^vx\
\prod_{m=1}^v\delta(x_m,\lambda_m)\nonumber\\[2.0ex] &
\le\ \int_{\cB_v(s\rho_1+(1-s)\rho_2)}\rd^vx\
\prod_{m=1}^v\delta(x_m,\lambda_m)\  =\ \alpha\left(s\rho_1+\left(1-s\right)\rho_2;\lambda\right)\,.
\end{align}
\end{proof}

\vskip 0.3cm

Now, log--concavity is a local property of $\alpha(\rho;\lambda)$,
yet it brings global information about the conditional moments of
$X$. To see this, we observe that since $\alpha(\rho;\lambda)$ is
twice differentiable with respect to~$\rho$,
eq.~(\ref{eq:logconcavity}) is equivalent to  
\begin{equation}
\alpha\partial_{\rho}^2\alpha - (\partial_\rho\alpha)^2 \le 0\,.
\label{eq:logconcequiv}
\end{equation}
We iterate eq.~(\ref{eq:recursion}) to express the above derivatives in
terms of conditional expectations. In first place, evaluating
that equation at $n=0$ yields
\begin{equation}
\partial_\rho \alpha = \frac{v}{2\rho}\alpha -
\frac{\alpha}{2\rho}\sum_{k=1}^v \frac{\E[X_k^2\,|\,X\in\cB_v(\rho)]}{\lambda_k}\,.
\label{eq:frstder}
\end{equation}
Property $p_1$) then implies
\begin{equation}
\sum_{k=1}^v \frac{\E[X_k^2\,|\,X\in\cB_v(\rho)]}{\lambda_k}\le v\,.
\label{eq:lowestineq}
\end{equation}
Though trivial, eq.~(\ref{eq:lowestineq}) calls for two remarks. The first one is that a
sufficient (but not necessary) condition for it to hold true is
$\E[X_k^2\,|\,X\in\cB_v(\rho)]\le \lambda_k$ \ $\forall k$,  which has already been established.  
In second place, differentiating $\alpha$ in $\rho$ an arbitrary number of
times generates always symmetric expressions with respect to the directional
indices, since $\rho$ is not tied to any specific direction. In
particular, this is the case with the second derivative,
\begin{equation}
\partial^2_\rho\alpha = \frac{\alpha}{\rho^2}\biggl\{\frac{v(v-2)}{4} -
\frac{v}{2}\sum_{k=1}^v\frac{\E[X_k^2\,|\,X\in\cB_v(\rho)]}{\lambda_k} +\frac{1}{4}\sum_{k,j=1}^v\frac{\E[X_j^2X_k^2\,|\,X\in\cB_v(\rho)]}{\lambda_j\lambda_k}\biggr\}\,.
\label{eq:scndder}
\end{equation}
We see that all directional indices are again summed over. We shall come back in sect.~4
to the rational coefficients multiplying the expectation values on the
\rhs of eqs.~(\ref{eq:frstder}) and (\ref{eq:scndder}). For the time being, we
finalize our argument by inserting these expressions into
eq.~(\ref{eq:logconcequiv}). A little algebra yields
\begin{equation}
\frac{\alpha^2}{4\rho^2}\left\{\sum_{k=1}^v \frac{\var\bigl(X_k^2\,|\,X\in\cB_v(\rho)\bigr)}{\lambda_k^2} -
  2v + \sum_{j\ne
    k}\frac{\cov\bigl(X^2_j,X^2_k\,|\,X\in\cB_v(\rho)\bigr)}{\lambda_j\lambda_k}
\right\}\le 0\,.
\label{eq:finallogconc}
\end{equation}
Eq.~(\ref{eq:finallogconc}) describes the log--concavity of $\alpha$ in terms of
conditional expectations. Did we not know about  Pr\'ekopa's
Theorem, we could regard it as a result of eqs.~(\ref{eq:varineq}) and
(\ref{eq:covineq}). Unfortunately, the converse does not hold: it is not possible to infer
eqs.~(\ref{eq:varineq}) and (\ref{eq:covineq}) from
eq.~(\ref{eq:finallogconc}), as contributions along different directions could
compensate while keeping the \lhs negative.  Nevertheless, if
eqs.~(\ref{eq:varineq}) and (\ref{eq:covineq}) were simultaneously violated
for all indices, $\alpha(\rho;\lambda)$ could not be logarithmic concave at
all. Therefore, eq.~(\ref{eq:finallogconc}) tells us that at least some of the
correlation  inequalities must hold. To conclude, a full proof of
eq.~(\ref{eq:varineq}) cannot follow from the property of log--concavity, so
we need to look elsewhere.  

\section{Variance reduction in the regime of strong truncation}

For the sake of conciseness, throughout this section we denote
conditional expectations by $\E[\,\cdot\,]$ instead of
$\E[\ \cdot\,\,|\,X\in\cB_v(\rho)]$. Our starting point consists in regarding
eq.~(\ref{eq:varineq}) as an upper bound to $\E[X_n^4]$. This suggests to
consider the wider inequality chain   
\begin{equation}
\E[X_n^4] \le \E[X_n^2]\left(2\lambda_n  + \E[X_n^2]\right) \le
\lambda_n\left(2\lambda_n  + \E[X_n^2]\right)\le 3\lambda_n^2\,.
\label{eq:ineqchain}
\end{equation}
The leftmost bound is in fact a recast of eq.~(\ref{eq:varineq}). If for a
moment we give it for granted, the second and third ones follow as an
immediate consequence of $\E[X_n^2]\le\lambda_n$. Although our final target is
just represented by eq.~(\ref{eq:varineq}), it makes sense to first consider
the two rightmost bounds: if they turn out to be violated,
eq.~(\ref{eq:varineq}) cannot be correct. The loosest one is once
more the trivial inequality $\E[X_n^4]\le 3\lambda_n^2$, which we have
already established. By contrast, the inequality    
\begin{equation}
\E[X_n^4] \le \lambda_n\left(2\lambda_n + \E[X_n^2]\right)\,
\label{eq:losineq}
\end{equation}
is less obvious. In sect. 3.1 we prove it. Our argument is
based on straightforward algebraic manipulations of the Gaussian integrals over $\cB_v(\rho)$.
We include it in the present note for a twofold reason: on the one hand it gives a
feeling of the optimality of eq.~(\ref{eq:varineq}), on the other it represents
the only general result we have, valid across the whole parameter space. 

\subsection{A loose yet general bound to $\E[X_n^4\,|\,X\in\cB_v(\rho)]$}

In order to prove eq.~(\ref{eq:losineq}), we use a standard trick, consisting
in a rescaling of $\lambda_n$ by an external parameter~$\tau$, so as to obtain the
moments of $X_n$ by differentiation of $\alpha$ in $\tau$. More precisely,
we introduce the function   
\begin{equation}
\cH(\tau) =
\frac{1}{\sqrt{\tau}}\,\alpha\left(\rho;\left\{\lambda_1,\dots,\frac{\lambda_n}{\tau},\dots,\lambda_v\right\}\right)
= \frac{1}{\sqrt{\tau}}\, \int_{\cB_v(\rho)}\rd^{v}x\ \delta(x_n;\lambda_n/\tau)\prod_{m\ne n}\delta(x_m,\lambda_m)\,,
\end{equation}
whose dependence upon $\rho$ and $\lambda$ we leave implicit. Differentiating $\cH(\tau)$ under the integral sign yields 
\begin{equation}
\E[X_n^{2k}] =
(-1)^k\frac{2^k\lambda_n^{k}}{\alpha}\frac{\partial^k\cH}{\partial\tau^k}\biggr|_{\tau=1}\,,\qquad k=0,1,2,\ldots
\label{eq:Xkmoms}
\end{equation}
At the same time, derivatives of $\cH(\tau)$ can be taken via the chain rule
of differentiation, which allows us to express them as algebraic combinations
of $\alpha$ and its derivatives in $\lambda_n$. For instance, with regard to
the second and fourth moments, we find
\begin{align}
\label{eq:Hderiv1}
\frac{\partial\cH}{\partial\tau} & = -\frac{1}{2\tau^{3/2}}\left(\alpha +
  2\lambda_n\partial_n\alpha\right)\,,\\[2.0ex]
\label{eq:Hderiv2}
\frac{\partial^2\cH}{\partial \tau^2} & = \frac{1}{4\tau^{5/2}}\left(3\alpha +
  8\lambda_n\partial_n\alpha + 4\lambda_n^2\partial_n^2\alpha\right)\,.
\end{align}
Consider first the lowest order derivative. By inserting eq.~(\ref{eq:Hderiv1}) into
eq.~(\ref{eq:Xkmoms}) evaluated at $k=1$, we obtain $\E[X_n^2] =
\lambda_n[1+2(\lambda_n/\alpha)\partial_n\alpha]$. Comparing this with
$\E[X_n^2] = \lambda_n(\alpha_n/\alpha)$ yields 
\begin{equation}
\alpha_n = \alpha + 2\lambda_n \partial_n\alpha\,.
\label{eq:alphakalpha}
\end{equation}
Eq.~(\ref{eq:alphakalpha}) coincides with eq.~(\ref{eq:derlk}) evaluated at
$n_1=\ldots=n_v=0$. Owing to property $p_2$) of sect.~2, we infer
$\alpha_n\le\alpha$ and thus we find again $\E[X_n^2] =
\lambda_n(\alpha_n/\alpha)\le\lambda_n$. Consider then the fourth moment. If
we insert 
eq.~(\ref{eq:Hderiv2}) into 
eq.~(\ref{eq:Xkmoms}) evaluated at $k=2$, and then make use of
eq.~(\ref{eq:alphakalpha}), we easily arrive at
\begin{equation}
\E[X_n^4] = 4\lambda_n\E[X_n^2] - \lambda_n^2 + 4\lambda_n^4\frac{\partial^2_n\alpha}{\alpha}\,.
\end{equation}
In order to estimate $\partial^2_n\alpha$, we differentiate
both sides of eq.~(\ref{eq:alphakalpha}) with respect to
$\lambda_n$. We then invoke again property $p_2$) of sect. 2, thus obtaining 
\begin{equation}
\partial^2_n\alpha = \frac{1}{2\lambda_n}\left(\partial_n\alpha_n -
3\partial_n\alpha\right)\le-\frac{3}{2\lambda_n}\partial_n\alpha = -\frac{3}{4\lambda_n^2}(\alpha_n-\alpha)
= -\frac{3\alpha}{4\lambda_n^3}\left\{\E\left[X_n^2\right] - \lambda_n\right\}\,.
\end{equation}
This estimate is sufficient to prove eq.~(\ref{eq:losineq}). 

\subsection{First argument in favor of eq.~(\ref{eq:varineq})}

\begin{center}
  \begin{figure}[!t]
    \begin{center}
      \begin{tabular}{p{0.45\textwidth}p{0.45\textwidth}}
        \hskip 0.5cm\includegraphics[width=0.40\textwidth]{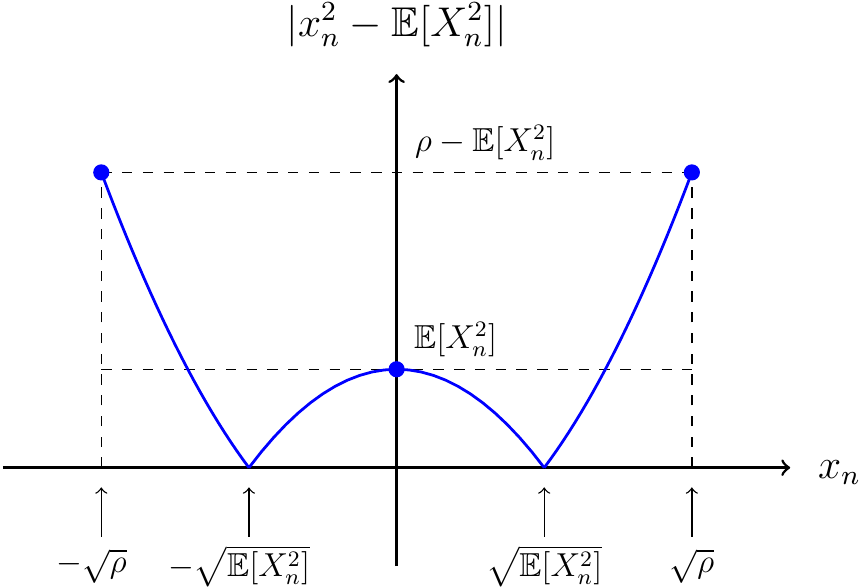}
        \begin{center}
          \vskip -0.3cm
          {\it a)}\hskip 0.3cm \phantom{ccc}
          \end{center}
        &
        \hskip 0.3cm\includegraphics[width=0.40\textwidth]{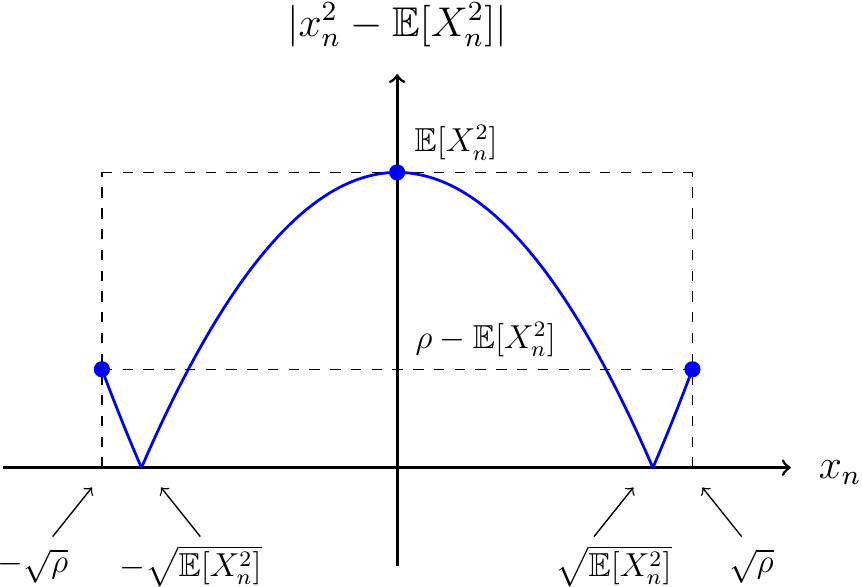}
        \begin{center}
          \vskip -0.3cm
          {\it b)}\hskip 0.3cm \phantom{cc}
          \end{center}
      \end{tabular}
      \vskip -0.5cm
      \caption{\small Qualitative behavior of the function $x_n \mapsto |x_n^2
            - \E[X_n^2]|$ within the support
            $(-\sqrt{\rho},+\sqrt{\rho})$. Note that, depending on  $\rho$ and 
            $\E[X_n^2]$, the function might have one or two
            maxima. Specifically,  plot~{\it a)} represents the
            function if $\rho-\E[X_n^2]>\E[X_n^2]$, while plot {\it
              b)} applies if $\rho-\E[X_n^2]<\E[X_n^2]$.\label{fig:esssup}} 
    \end{center}
  \vskip -0.3cm
  \end{figure}
\end{center}

In the regime of strong truncation, eq.~(\ref{eq:varineq}) can be inferred from
H\"older's inequality. We recall that if $p,q>1$ are two numbers satisfying
$1/p+1/q=1$ and $X$, $Y$ are stochastic variables on a
given probability space, then $\E[|XY|]\le(\E[|X|^p])^{1/p}(\E[|Y|^q])^{1/q}$.
In our case, we have 
\begin{align}
\var(X_n^2) & = \E\left[\left(X_n^2-\E[X_n^2]\right)^2\right] =
\E\left[X_n^4\right] - \E\left[X_n^2\right]^2 = \E\left[X_n^4 -
  \E[X_n^2]^2\right] \nonumber\\[2.0ex]
& = \E\left[\left(X_n^2 - \E[X_n^2]\right)\left(X_n^2 +
    \E[X_n^2]\right)\right] \le \E\left[\left|\left(X_n^2 - \E[X_n^2]\right)\left(X_n^2 +
    \E[X_n^2]\right)\right|\right] \nonumber\\[2.0ex]
& \le \left\{\E\left[\left|X_n^2 - \E[X_n^2]\right|^p\right]\right\}^{1/p} \left\{\E\left[\left(X_n^2 + \E[X_n^2]\right)^q\right]\right\}^{1/q}\,.
\end{align}
The latter inequality holds true for any finite choice of $p,q$, provided their
reciprocals sum to one. Accordingly, it holds as well in the joint limit $q\to
1^+$, $p= q/(q-1)\to\infty$, where it amounts to
\begin{equation}
\var(X_n^2) \le 2\,h\,\E[X_n^2]\,,
\end{equation}
with
\begin{equation}
h \equiv \lim_{p\to\infty} \left\{\E\left[\left|X_n^2 - \E[X_n^2]\right|^p\right]\right\}^{1/p} =
\,\text{ess\,sup} \left(\left|X_n^2 - \E[X_n^2]\right|\right)\,.
\label{eq:hdef}
\end{equation}
Recall that the essential supremum of a real--valued function $f$ is defined
by ${\rm ess\,sup}\, f \equiv \inf\{a\in\RR:\ \mu(\{x:f(x)>a\})=0\}$.  In
particular, the measure $\mu$ which is understood in eq.~(\ref{eq:hdef}) is
the marginal probability measure of $X_n$, 
{\it i.e.}
\begin{align}
\rd \mu(x_n) =
\frac{\alpha^{(v-1)}(\rho-x_n^2;\lambda_{(n)})\,\delta(x_n,\lambda_n)}{\alpha^{(v)}(\rho;\lambda)}\ \rd x_n\,.
\label{eq:margprob}
\end{align}
Owing to the modulating factor $\alpha^{(v-1)}(\rho-x_n^2;\lambda_{(n)})$,
$\mu$ is neither Gaussian nor log--concave (in $x_n$). In sect.~4 we shall say more
about eq.~(\ref{eq:margprob}) and the factorization of its numerator
into functions of resp. $\lambda_n$ and $\lambda_{(n)}$. For the time being, we observe that $\mu$
has support in the interval $(-\sqrt{\rho},+\sqrt{\rho})$. Depending on how $\E[X_n^2]$
compares with $\rho/2$, $h$ might assume one of the values $h_1 = \rho -
\E[X_n^2]$ or $h_2 = \E[X_n^2]$, as qualitatively represented in
Fig.~\ref{fig:esssup}. As far as we are concerned, we do not need to
establish which among Fig.~\ref{fig:esssup}a and Fig.~\ref{fig:esssup}b
provides the correct qualitative behavior for $x_n\mapsto\left|x_n^2 -
\E[X_n^2]\right|$: numerical computations suggest that Fig.~\ref{fig:esssup}b
is not realized for any choice of $n$, $\rho$ and $\lambda$, yet this
information is irrelevant for what follows. More precisely, we distinguish three cases: 
\begin{itemize}
\item[{\it i)}]{$\rho\le\lambda_n$ ({\it strong truncation}): in this case $h\le\lambda_n$. Indeed, since
    $\E[X_n^2]\le\lambda_n$, both $h_1$ and $h_2$ lie
    below $\lambda_n$. In this region, we have no analytic argument in favour of
    Fig.~\ref{fig:esssup}a or Fig.~\ref{fig:esssup}b.} 
\item[{\it ii)}]{$\rho> 2\lambda_n$ ({\it weak truncation}): in this case
    $h>\lambda_n$. Indeed, again from $\E[X_n^2]\le\lambda_n$, we deduce
    $\rho-\E[X_n^2]\ge\rho-\lambda_n>\lambda_n\ge\E[X_n^2]$. Here, the correct
    profile of $\left|x_n^2-\E[X_n^2]\right|$ is certainly the one depicted in Fig.~\ref{fig:esssup}a.}
\item[{\it iii)}]{$\lambda_n<\rho<2\lambda_n$: in this case it is difficult to
    conclude anything about $h$, except that by continuity there exists a
    value $\lambda_n<\rho_*(\lambda_{(n)})<2\lambda_n$, possibly depending on
    $\lambda_{(n)}$, such that $h\le\lambda_n$ $\Leftrightarrow$ 
    $\rho\le\rho_*$.}
\end{itemize}

To conclude, the estimate obtained from H\"older's inequality
is certainly as strict as needed for eq.~(\ref{eq:varineq}) to hold true only
in case of strong truncation, {\it i.e.} for $\rho\le\lambda_n$. In addition, there is a
crossover region where the same estimate might be sufficiently strict, while it
becomes definitely too loose in the region of weak truncation.

\subsection{Second argument in favour of eq.~(\ref{eq:varineq})}

In order to extend the above proof to the region
$0<\rho\le 2\lambda_n$, we work on the integral
representations of $\E[X_n^4] = \lambda_n^2(\alpha_{nn}/\alpha)$ and $\E[X_n^2]
=\lambda_n(\alpha_{n}/\alpha)$. In terms of these, $\Delta_n$ reads
\begin{equation}
\Delta_n = \frac{\lambda_n^2}{\rho^2}\left[\frac{\alpha_{nn}}{\alpha} -
  \left(\frac{\alpha_n}{\alpha}\right)^2 - 2 \frac{\alpha_n}{\alpha}\right]\,.
\label{eq:Deltak}
\end{equation}
Now, we observe that independently of $v$, $\alpha_{n:k}$ is bounded from
above by $(\rho/\lambda_n)^{k-p}\alpha_{n:p}$ for any $p<k$. Indeed, since
$x\in\cB_v(\rho)\ \Rightarrow\ -\sqrt{\rho}<x_n<\sqrt{\rho}$, we have
\begin{equation}
\alpha_{n:k} = \frac{\rho^k}{\lambda_n^k}\int_{\cB_v(\rho)}\rd^vx\
\frac{x_n^{2k}}{\rho^k}\prod_{m=1}^v\delta(x_m,\lambda_m) \le 
\frac{\rho^{k}}{\lambda_n^k} \int_{\cB_v(\rho)}\rd^v x\
\frac{x_n^{2p}}{\rho^p}\prod_{m=1}^v\delta(x_m,\lambda_m) = \frac{\rho^{k-p}}{\lambda_n^{k-p}}\alpha_{n:p}\,.
\label{eq:lowdominance}
\end{equation}
Thus, we immediately obtain
\begin{equation}
\Delta_n \le \frac{\lambda_n^2}{\rho^2}\left[\left(\frac{\rho}{\lambda_n} - 2\right)\frac{\alpha_n}{\alpha} -
\left(\frac{\alpha_n}{\alpha}\right)^2\right] \le 0\,,\qquad \text{if }\ \rho \le 2\lambda_n\,.
\end{equation}
This conclusion is somewhat conservative, as indeed $\Delta_n\le 0\
\Leftrightarrow\ \rho \le \rho_*$, being $\rho_*$ implicitly defined by the
non--linear equation $\rho_* =  2\{\lambda_n +
\E[X_n^2\,|\,X\in\cB_v(\rho_*)]\}$. By continuity, the latter is certainly
fulfilled by some $2\lambda_n<\rho_*\le 4\lambda_n$. The argument
presented here does not apply for $\rho>4\lambda_n$.

\section{Weak truncation expansion}

In order to study the variance of the square components of $X$ in the
regime of weak truncation, we need to develop an appropriate
formalism. To start with, we observe that the constraint
$X\in\cB_v(\rho)$ becomes increasingly unrestrictive as
$\rho\to\infty$. As a consequence, we have the asymptotic
factorization  
\begin{equation}
\alpha^{(v)}_{1:k_1\, 
  \ldots\, v:k_v}(\rho;\lambda) \ \
\widesim{\ \rho\,\gg\,\max_j\{\lambda_j\}} \ \ \
\prod_{j=1}^v\alpha^{(1)}_{j:k_j}(\rho;\lambda_{j})\,.
\label{eq:asymp1}
\end{equation}
The larger is $\rho$, the less is the error made in approximating
$\alpha_{1:k_1\ldots v:k_v}$ by its factorized counterpart. We aim at
characterizing the corrections to eq.~(\ref{eq:asymp1}) when  $\rho$ is
large yet finite. Actually, we are not interested in a complete 
factorization of the Gaussian integrals: if $\rho\gg\lambda_n$ just for some
$1\le n\le v$, we look at the partial factorization occurring along the
$n^{\rm th}$ direction. Note that: {\it i})~in the regime of weak
truncation, every rational combination of Gaussian integrals --- such as
$\Delta_n$ --- is led by its factorized counterpart; as we shall see,
the latter is subject to relevant simplifications in case of ratios of integrals;
{\it ii}) 1--dimensional integrals cannot be further simplified, as
they amount to lower incomplete gamma functions, 
\begin{align}
\alpha^{(1)}_{n:k}(\rho;\lambda_n) =
\frac{\,2^k\,}{\sqrt{\pi}}\,\gamma\left(k+\frac{1}{2},\frac{\rho}{2\lambda_n}\right)\,,\qquad\qquad
\gamma(s,x) = \int_0^x\rd t\, t^{s-1}\re^{-t}\,.
\label{eq:gammafunc}
\end{align}

\subsection{Expansion of Gaussian integrals}

\begin{center}
\begin{figure}
\begin{center}
\psfrag{tag1}{\footnotesize $(v-1)$--dimensional}
\psfrag{tag2}{\footnotesize hyperplane $\perp \ y_i$}
\psfrag{tag3}{\small $y_i$}
\includegraphics[angle=-22,width=2.6in]{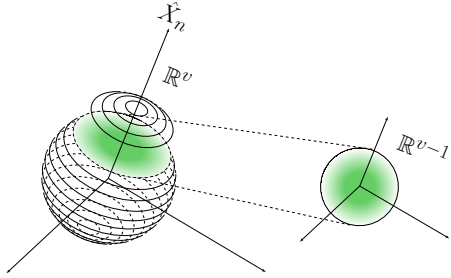}
\vskip -1.6cm
\caption{\small Slicing $\cB_v(\rho)$ orthogonally to the $n^{\rm th}$
  coordinate axis $\hat X_n$ amounts to representing it as an uncountable union of 
  Euclidean balls, living in the $(v-1)$--dimensional subspace orthogonal to
  $\hat X_n$, with square radius varying from $0$ to $\rho$.\label{fig:slicing}}    
\end{center}
\end{figure}
\end{center}

In order to present the idea, we first focus
 on $\alpha$. If $\rho\gg\lambda_n$ for some $1\le n\le v$, we
slice the integration domain orthogonally to the $n^{\rm th}$
direction, as depicted in Fig.~\ref{fig:slicing}. From a geometrical
point of view, this corresponds to representing  $\cB_v(\rho)$ as an
uncountable union of $(v-1)$--dimensional Euclidean balls, {\it i.e.}
\begin{align}
\cB_v(\rho) = \!\!\!\!\!
  \bigcup_{x_n\in (-\sqrt{\rho},+\sqrt{\rho})}\!\!\!\!\! \left\{y\in \RR^v:\ y_n=x_n\,,\ \ y_{(n)}\in
  \cB_{v-1}(\rho - x_n^2)\right\}\,.
\label{eq:slicing}
\end{align}
Such technique has been first considered by Ruben~\cite{ruben} with
the aim of obtaining an integral recurrence relationship 
on the dimensionality of $\alpha$. Interpreting the integration
domain in terms of eq.~(\ref{eq:slicing}) indeed yields
\vskip -0.5cm
\begin{equation}
\alpha^{(v)}(\rho;\lambda) = \int_{-\sqrt{\rho}}^{+\sqrt{\rho}} \rd x_n\
\delta(x_n,\lambda_n)\ \alpha^{(v-1)}(\rho-x_n^2;\lambda_{(n)})\,.
\label{eq:alphasliced}
\end{equation}
Since $\alpha(\rho-x_n^2;\lambda_{(n)})$ is a smooth function of its
first argument $\rho-x_n^2$, we propose to expand it in Taylor series around $x_n^2=0^+$,
\begin{align}
\alpha^{(v-1)}(\rho-x_n^2;\lambda_{(n)})
& = \sum_{k=0}^{\infty}\frac{(-1)^k}{k!}\lambda_n^{k}\left(\frac{x_n^2}{\lambda_n}\right)^k\partial^k_\rho
\alpha^{(v-1)}(\rho;\lambda_{(n)})\nonumber\\[2.0ex]
& = \alpha^{(v-1)}(\rho;\lambda_{(n)})\sum_{k=0}^{\infty}\frac{(-1)^k}{k!}\left(\frac{\lambda_n}{\rho}\right)^k\left(\frac{x_n^2}{\lambda_n}\right)^k\eta^{(v-1)}_k(\rho;\lambda_{(k)})\,,
\label{eq:alphataylor}
\end{align}
with the functions $\eta_k$ defined by
\begin{equation}
\eta^{(v)}_k(\rho;\lambda) = \left\{\begin{array}{ll}1\,, & \qquad k=0\,,\\[2.0ex]
\displaystyle{[{\alpha^{(v)}(\rho;\lambda)}]^{-1}\rho^k\partial_\rho^k\alpha^{(v)}(\rho;\lambda)}\,,&\qquad
k\ge 1\,.\end{array}\right.
\end{equation}
When inserted into eq.~(\ref{eq:alphasliced}), the Taylor series turns into a
weak truncation expansion of $\alpha$, namely
\begin{align}
\alpha^{(v)}(\rho;\lambda)\, & =\, \alpha^{(v-1)}(\rho;\lambda_{(n)}) \sum_{k=0}^\infty
\frac{(-1)^k}{k!}\,\left(\frac{\lambda_n}{\rho}\right)^k\,
\alpha_{n:k}^{(1)}(\rho;\lambda_n)\,\eta_k^{(v-1)}(\rho;\lambda_{(n)})
\nonumber\\[1.5ex] 
&=  \alpha^{(1)}(\rho;\lambda_n)\,\alpha^{(v-1)}(\rho;\lambda_{(n)})
\, -\, \frac{\lambda_n}{\rho}\, \alpha_n^{(1)}(\rho;\lambda_n)\,\alpha^{(v-1)}(\rho;\lambda_{(n)})\,\eta_1^{(v-1)}(\rho;\lambda_{(n)})\nonumber\\[2.0ex]
& +\,
\frac{1}{2}\,\frac{\lambda_n^2}{\rho^2}\,\alpha_{nn}^{(1)}(\rho;\lambda_n)\,\alpha^{(v-1)}(\rho;\lambda_{(n)})\eta_2^{(v-1)}(\rho;\lambda_{(n)})
\, +\, {\rm O}\left(\frac{\lambda_n^3}{\rho^3}\right)\,.
\label{eq:weaktrunc1}
\end{align}
Although a complete factorization into functions of $\lambda_n$ and
$\lambda_{(n)}$ is not exactly realized at finite $\rho$, we see that
it occurs at each order of the expansion. We warn that eq.~(\ref{eq:weaktrunc1})
has been obtained upon bringing an infinite sum under an integral
sign. Such exchange of limits is delicate, so it is not {\it a priori}
evident whether the resulting expansion converges or approximates its
target just as an asymptotic series. We shall come back to this point later
on. We also stress that, while power counting is performed by
factors of $(\lambda_n/\rho)^k$, additional powers and exponentially
small terms in $\rho$ are still hidden within the coefficient
functions\footnote{The reader will observe that the coefficient
  function $\alpha^{(v-1)}$ showing up in each term of the expansion
  is totally useless, as it simplifies with the one attached to
  $\eta_k^{(v-1)}$. Such redundancy is real, yet it turns useful when
  ratios of Gaussian integrals are considered, as we shall see in
  sects.~4.2 and 5.} $\alpha_{n:k}^{(1)}$ and $\eta_{k}^{(v-1)}$.    

To simplify the notation, in the sequel we drop all function arguments,
whenever this does not generate confusion. Thus, we shorten
eq.~(\ref{eq:weaktrunc1}) to
\begin{equation}
\alpha^{(v)} =\, \alpha^{(1)}\,\alpha^{(v-1)}\, -\, \frac{\lambda_n}{\rho}\, \alpha_n^{(1)}\alpha^{(v-1)}\,\eta^{(v-1)}_1\,+\,
\frac{1}{2}\frac{\lambda_n^2}{\rho^2}\,\alpha_{nn}^{(1)}\,\alpha^{(v-1)}\,\eta^{(v-1)}_2 \, +\, {\rm O}\left(\frac{\lambda_n^3}{\rho^3}\right)\,.
\label{eq:a0exp}
\end{equation}
The same technique can be straightforwardly applied to
$\alpha_{n:p}$. For instance, we have for $p=1,2,\ldots$ 
\begin{align}
\alpha_{n}^{(v)} & =\, \alpha_{n:1}^{(1)}\,\alpha^{(v-1)}\, -\, \frac{\lambda_n}{\rho}\,\alpha_{n:2}^{(1)}\,\alpha^{(v-1)}\,\eta^{(v-1)}_1\, +\,
\frac{1}{2}\,\frac{\lambda_n^2}{\rho^2}\,\alpha_{n:3}^{(1)}\alpha^{(v-1)}\,\eta^{(v-1)}_2 \,
+\, {\rm O}\left(\frac{\lambda_n^3}{\rho^3}\right)\,,\label{eq:a1exp}\\[1.0ex]
\alpha_{nn}^{(v)} & =\, \alpha_{n:2}^{(1)}\,\alpha^{(v-1)}\, -\, \frac{\lambda_n}{\rho}\,\alpha_{n:3}^{(1)}\,\alpha^{(v-1)}\,\eta^{(v-1)}_1\, +\,
\frac{1}{2}\,\frac{\lambda_n^2}{\rho^2}\,\alpha_{n:4}^{(1)}\alpha^{(v-1)}\,\eta^{(v-1)}_2
\,+\, {\rm
  O}\left(\frac{\lambda_n^3}{\rho^3}\right)\,,\label{eq:a2exp}\\[1.0ex]
& \hskip 0.15cm \vdots\nonumber
\end{align}
Since the above expansions are all based on
eq.~(\ref{eq:alphataylor}), the coefficient functions $\eta_k$ are 
the same independently of $p$. By contrast, the multiplicity of
the index $n$ of the 1--dimensional integrals contributing to each order is
shifted forward as $p$ increases. 

Now, when it comes to expanding Gaussian integrals with more than one index,
the above procedure is carried out in a slightly different way. For
instance, in order to expand $\alpha_{nm}$ we need to take into account
factors of $x_n^2$ and $x_m^2$ under the integral sign. Accordingly, we slice $\cB_v(\rho)$
subsequently along the $n^{\rm th}$ and $m^{\rm th}$ directions under the
assumption $\rho\gg \max\{\lambda_n,\lambda_m\}$. The analogous of eq.~(\ref{eq:alphasliced}) reads
\begin{equation}
\alpha_{nm}^{(v)}(\rho;\lambda) = \int_{-\sqrt{\rho}}^{+\sqrt{\rho}} \rd x_n\
\frac{x_n^2}{\lambda_n}\delta(x_n,\lambda_n)\int_{-\sqrt{\rho}}^{+\sqrt{\rho}} \rd x_m\
\frac{x_m^2}{\lambda_m}\delta(x_m,\lambda_m)\ \alpha^{(v-2)}(\rho-x_n^2-x_m^2;\lambda_{(nm)})\,,
\label{eq:alphasliced2}
\end{equation}
with $\lambda_{(nm)}\equiv\{\lambda_k\}_{k\ne n,m}$. Again, we expand
$\alpha(\rho-x_n^2-x_m^2;\lambda_{(nm)})$ in Taylor series around the point
$x_n^2+x_m^2=0^+$, thus obtaining 
\begin{align}
\alpha^{(v-2)}(\rho-x^2_n-x^2_m;\lambda_{(nm)}) \ & = \ \sum_{j=0}^\infty
\frac{(-1)^j}{j!}(x^2_n+x^2_m)^j \partial^j_\rho\alpha^{(v-2)}(\rho;\lambda_{(nm)})
\nonumber \\[0.0ex]
& \hskip -3.0cm =\  \sum_{j=0}^\infty
\frac{(-1)^j}{j!}\sum_{k=0}^j\binom{j}{k}x^{2k}_nx^{2(j-k)}_m \partial^j_\rho\alpha^{(v-2)}(\rho;\lambda_{(nm)}) \nonumber \\[0.0ex]
& \hskip -3.0cm =\  \alpha^{(v-2)}(\rho;\lambda_{(nm)}) \sum_{j=0}^\infty
\frac{(-1)^j}{j!}\sum_{k=0}^j\binom{j}{k}\left(\frac{\lambda_n}{\rho}\right)^k\left(\frac{\lambda_m}{\rho}\right)^{j-k}\nonumber\\[0.0ex]
& \hskip 2.3cm \cdot \biggl(\frac{x^{2}_n}{\lambda_n}\biggr)^k\left(\frac{x^{2}_m}{\lambda_m}\right)^{j-k}
\eta_j^{(v-2)}(\rho;\lambda_{(nm)})\,.
\label{eq:alphataylor2}
\end{align}
Here we have also used Newton's binomial formula to express
each term of the series in products of powers of $x_n^2$ and
$x_m^2$. Inserting this expression back into eq.~(\ref{eq:alphasliced2})
yields
\begin{equation}
\alpha_{nm}^{(v)} = \alpha_n^{(1)}\,\alpha_m^{(1)}\alpha^{(v-2)} - \frac{\lambda_n}{\rho}\,\alpha_{n:2}^{(1)}\,\alpha_m^{(1)}\,\alpha^{(v-2)}\,\eta_1^{(v-2)}
-\frac{\lambda_m}{\rho}\,\alpha_n^{(1)}\,\alpha_{m:2}^{(1)}\,\alpha^{(v-2)} \,\eta^{(v-2)}_1
+ {\rm O}\left(\frac{\lambda^2}{\rho^2}\right)\,.
\label{eq:a2exp2}
\end{equation}
Eq.~(\ref{eq:alphataylor2}) can be also used 
to obtain alternative expansions of $\alpha_{n:p}$ if
$\rho\gg\max\{\lambda_n,\lambda_m\}$ for some $m$, {\it e.g.}
\begin{align}
\alpha^{(v)} &\, =\, \alpha^{(1)}\,\alpha^{(1)}\alpha^{(v-2)} -
\frac{\lambda_n}{\rho}\,\alpha_{n:1}^{(1)}\,\alpha^{(1)}\,\alpha^{(v-2)}\,\eta^{(v-2)}_1 -
\frac{\lambda_m}{\rho}\,\alpha_{m:1}^{(1)}\,\alpha^{(1)}\,\alpha^{(v-2)}\,\eta^{(v-2)}_1 +
{\rm O}\left(\frac{\lambda^2}{\rho^2}\right)\,, 
\label{eq:a0exp2}
\end{align}
\vskip -0.4cm
\begin{align}
\alpha_{n}^{(v)} &\, =\, \alpha_{n:1}^{(1)}\,\alpha^{(1)}\,\alpha^{(v-2)} -
\frac{\lambda_n}{\rho}\,\alpha_{n:2}^{(1)}\,\alpha^{(1)}\,\alpha^{(v-2)}\,\eta^{(v-2)}_1 -
\frac{\lambda_m}{\rho}\,\alpha_{n:1}^{(1)}\,\alpha_{m:1}^{(1)}\,\alpha^{(v-2)}\,\eta^{(v-2)}_1
+ {\rm O}\left(\frac{\lambda^2}{\rho^2}\right)\,,\label{eq:a1exp2}\\[0.0ex]
& \ \ \vdots\nonumber
\end{align}
\vskip -0.3cm

\subsection{Exercise: relative amplitude of variances and covariances}

We make use of the above expansions to qualitatively
compare the correlations $X_n^2$ has with itself and the other square
components of $X$ in the regime of weak truncation. This is rather
instructive, because it shows how analytic cancellations occur in the
proposed formalism. In addition, the exercise inspires the following
unproved 
\begin{conj}
If $X\sim\cN_{v}(0,\Lambda)$ with $\Lambda = \diag(\lambda)$, the
covariance matrix of the vector $\{X_n^2\}_{n=1}^v$ conditioned to
$\cB_v(\rho)$ is diagonally dominant, {\it i.e.}
\begin{equation}
\var(X_n^2\,|\,X\in\cB_v(\rho)) \ge \sum_{m\ne n}
\left|\cov(X_n^2,X_m^2\,|\,X\in\cB_v(\rho))\right|\,,\qquad \rho\in\RR_+\,.
\label{eq:conjvar}
\end{equation}
\end{conj}
\vskip-0.6cm\hfill\qed

\noindent Eq.~(\ref{eq:conjvar}) is supported with no exceptions by
extensive numerical tests. As observed in ref.~\cite{palombi2}, it
entails the inequality $v^{-1}\sum_{k=1}^v\mu_k\le {\rho}/{(v+2)}$,
where $\mu=\{\mu_k\}_{k=1}^v$ denotes the eigenvalue spectrum of the
covariance matrix ${\frak S}_\cB$.  

\begin{center}
  \begin{figure}[!t]
    \begin{center}
    \includegraphics[width=0.55\textwidth]{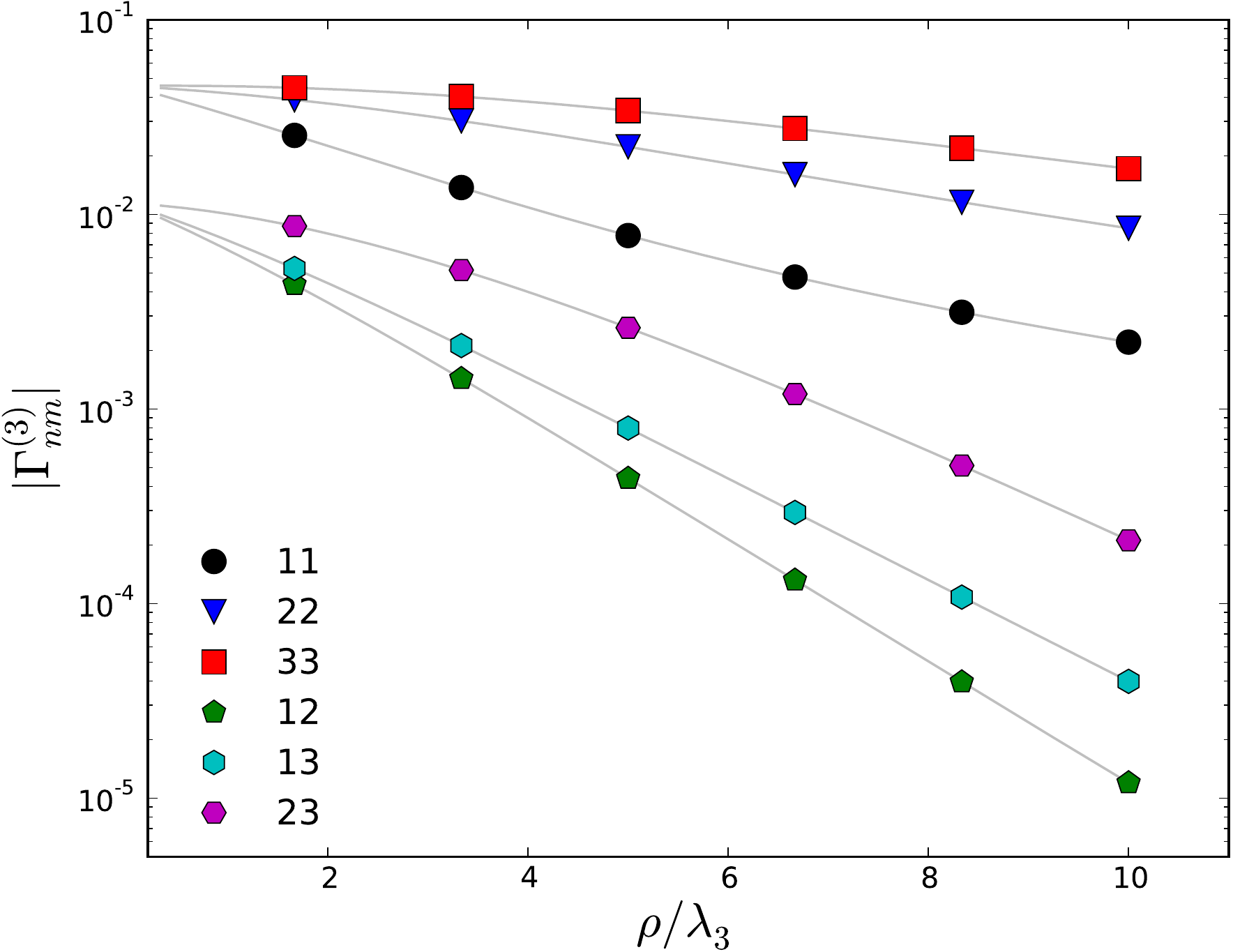}
    \caption{\small A plot of $|\Gamma^{(3)}_{nm}|$ vs. $\rho/\lambda_3$ at
          $\{\lambda_1,\lambda_2,\lambda_3\} = \{1,2,3\}$. Solid curves
          correspond to computed values: markers have been overimposed
          just to make the legend clear.\label{fig:Gammas}} 
    \end{center}
  \end{figure}
\end{center}

In terms of Gaussian integrals the scale invariant observables we focus on are
\begin{align}
\Gamma_{nn}^{(v)} & \, \equiv\, \frac{1}{\rho^2}\var(X_n^2\,|\,X\in\cB_v(\rho)) \, = \, \frac{\lambda_n^2}{\rho^2}\left[\frac{\alpha^{(v)}_{nn}}{\alpha^{(v)}} -
  \left(\frac{\alpha^{(v)}_n}{\alpha^{(v)}}\right)^2\right]\,,\label{eq:gammaii}\\[1.0ex]
\Gamma_{nm}^{(v)} & \, \equiv\, \frac{1}{\rho^2}\cov(X_n^2,X_m^2\,|\,X\in\cB_v(\rho)) \, = \,
\frac{\lambda_n}{\rho}\frac{\lambda_m}{\rho}\left[\frac{\alpha^{(v)}_{nm}}{\alpha^{(v)}}
  - \frac{\alpha^{(v)}_n}{\alpha^{(v)}}\frac{\alpha^{(v)}_m}{\alpha^{(v)}}\right]\,,\qquad
n\ne m\,.\label{eq:gammaik}
\end{align}
For illustrative purposes, we show in Fig.~\ref{fig:Gammas} a plot of
$|\Gamma^{(v)}_{nm}|$ vs. $\rho/\lambda_3$ at $v=3$ corresponding to
the choice $\{\lambda_1,\lambda_2,\lambda_3\} = \{1,2,3\}$. Both
$\Gamma_{nn}$ and $\Gamma_{nm}$ vanish as $\rho\to\infty$, yet the
former vanishes as $1/\rho^2$ due to the chosen normalization, whereas
the latter is exponentially damped. 
\vskip -0.45cm

\begin{figure}[!t]
  \centering
  \includegraphics[width=0.60\textwidth]{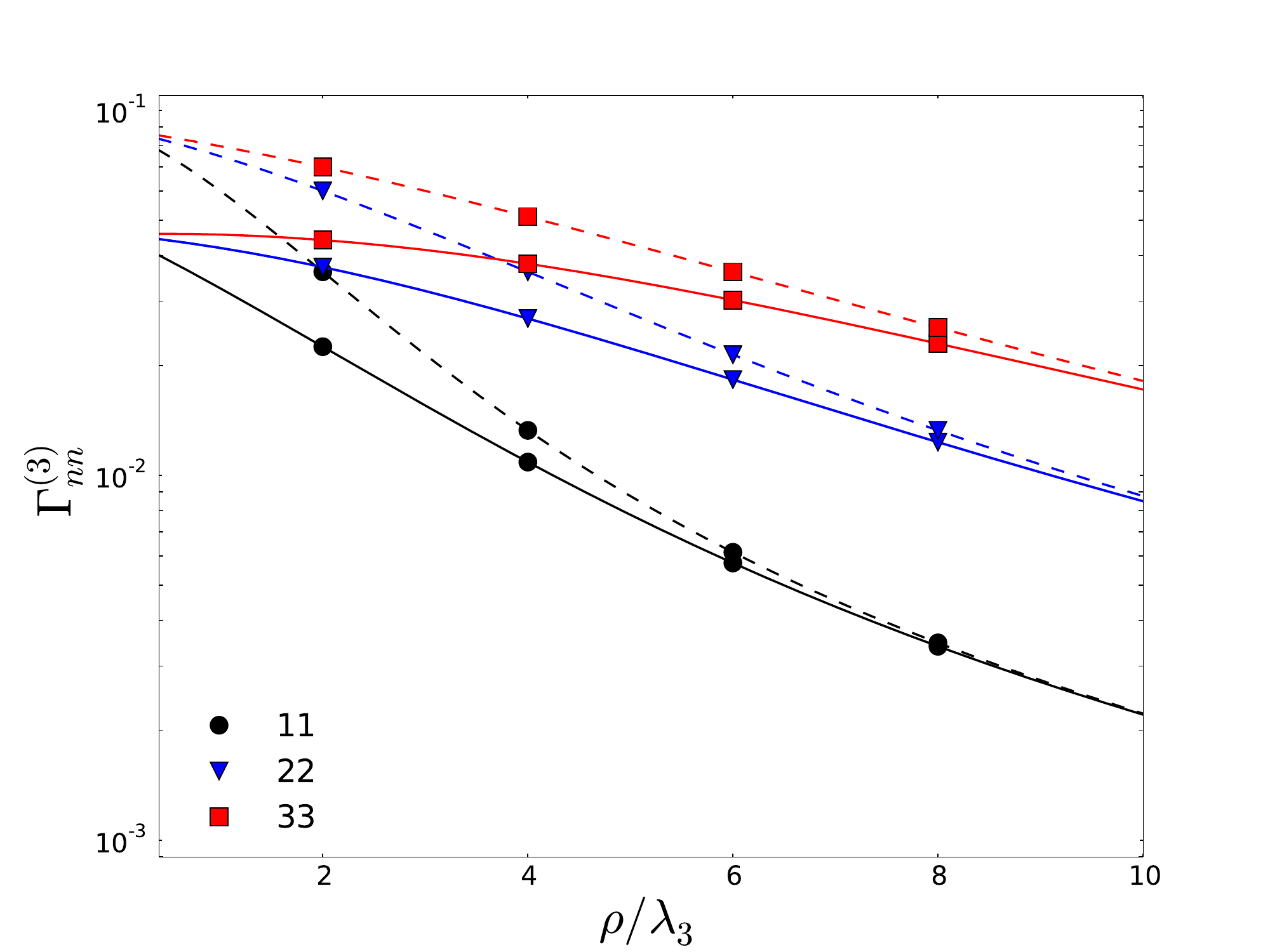}
  \caption{\small Convergence rate of $\Gamma^{(3)}_{nn}$ (solid lines)
    to $\Gamma^{(1)}_{nn}$ (dashed lines) at
    $\{\lambda_1,\lambda_2,\lambda_3\}=\{1,2,3\}$.
    Curves correspond to computed values: markers have been overimposed 
    just to make the legend clear.}
  \label{fig:conv} 
\end{figure}

We use eqs.~(\ref{eq:a0exp})--(\ref{eq:a2exp}) to work out the
expansion of $\Gamma_{nn}$ and
eqs.~(\ref{eq:a2exp2})--(\ref{eq:a1exp2}) for $\Gamma_{nm}$. In both
cases, in order to expand $\alpha^{-1}$ we rely on the Taylor formula
$(1-x)^{-1} = 1 + x + x^2 + {\rm O}(x^3)$. Thus, with regard to
$\Gamma_{nn}$ we have  
\begin{equation}
\frac{\alpha_{nn\phantom{:1}}^{(v)}}{\alpha^{(v)}} =
\frac{\alpha_{nn\phantom{:1}}^{(1)}}{\alpha^{(1)}} -
\frac{\lambda_n}{\rho}\left[\frac{\alpha_{n:3}^{(1)}}{\alpha^{(1)}} -
\frac{\alpha_{n:2}^{(1)}}{\alpha^{(1)}}\frac{\alpha_{n\phantom{:1}}^{(1)}}{\alpha^{(1)}}\right]\eta_1^{(v-1)}
+ {\rm O}\left(\frac{\lambda_n^2}{\rho^2}\right)\,,
\end{equation}
\begin{equation}
\left(\frac{\alpha_{n\phantom{:1}}^{(v)}}{\alpha^{(v)}}\right)^2 =
\left(\frac{\alpha_{n\phantom{:1}}^{(1)}}{\alpha^{(1)}}\right)^2 -
2\frac{\lambda_n}{\rho}\left[\frac{\alpha_{n:2}^{(1)}}{\alpha^{(1)}}\frac{\alpha_{n\phantom{:1}}^{(1)}}{\alpha^{(1)}}-\left(\frac{\alpha_{n\phantom{:1}}^{(1)}}{\alpha^{(1)}}\right)^3\right]\eta_1^{(v-1)} + {\rm O}\left(\frac{\lambda_n^2}{\rho^2}\right)\,,
\end{equation}
whence we obtain
\begin{equation}
\Gamma^{(v)}_{nn} = \Gamma^{(1)}_{nn} -\frac{\lambda_n^3}{\rho^3}\left[\frac{\alpha_{n:3}^{(1)}}{\alpha^{(1)}} -3
\frac{\alpha_{n:2}^{(1)}}{\alpha^{(1)}}\frac{\alpha_{n\phantom{:1}}^{(1)}}{\alpha^{(1)}} +2\left(\frac{\alpha_{n\phantom{:1}}^{(1)}}{\alpha^{(1)}}\right)^3\right]\eta_1^{(v-1)}
+ {\rm O}\left(\frac{\lambda_n^2}{\rho^2}\right)\,.  
\label{eq:weakvar}
\end{equation}
We see that the leading term of $\Gamma_{nn}$ coincides with its
1--dimensional counterpart. In particular, in Fig.~\ref{fig:conv} we show the 
rate at which $\Gamma^{(v)}_{nn}$ approaches $\Gamma^{(1)}_{nn}$ at $v=3$ and
 $\lambda = \{1,2,3\}$. Since
$\lim_{\rho\to\infty}\alpha_{n:k}^{(1)} = (2k-1)!!$, we have  
$\lim_{\rho\to\infty}\,(\rho^2/\lambda_n^2)\,\Gamma_{nn}^{(1)} = 2$,
and thus we find again
\begin{equation}
\Gamma^{(v)}_{nn}\ \ \widesim{\ \rho\,\gg\,\lambda_n} \ \ \
\frac{\lambda_n^2}{\rho^2}\left\{ 2  +  {\rm O}\left(\frac{\lambda_n}{\rho}\right)\right\}\,,
\label{eq:Gammaiiexp}
\end{equation}
apart from exponentially small terms in $\rho$. Analogously, we have
\begin{align}
& \frac{\alpha_{nm\phantom{:1}}^{(v)}}{\alpha^{(v)}} =
\frac{\alpha_{n\phantom{:1}}^{(1)}}{\alpha^{(1)}}\frac{\alpha_{m\phantom{:1}}^{(1)}}{\alpha^{(1)}} -
\frac{\lambda_n}{\rho}\left[\frac{\alpha_{n:2}^{(1)}}{\alpha^{(1)}}-\frac{\alpha_{m\phantom{:1}}^{(1)}}{\alpha^{(1)}}\left(\frac{\alpha_{n\phantom{:1}}^{(1)}}{\alpha^{(1)}}\right)^2\right]\,\eta_1^{(v-2)}\label{eq:discterm1}\nonumber \nonumber\\[0.0ex]
& \hskip 3.12cm -
\frac{\lambda_m}{\rho}\left[\frac{\alpha_{m:2}^{(1)}}{\alpha^{(1)}}-\frac{\alpha_{n\phantom{:1}}^{(1)}}{\alpha^{(1)}}\left(\frac{\alpha_{m\phantom{:1}}^{(1)}}{\alpha^{(1)}}\right)^2\right]\,\eta_1^{(v-2)}
+ {\rm O}\left(\frac{\lambda^2}{\rho^2}\right)\,,
\end{align}
\begin{align}
& \frac{\alpha_{n\phantom{:1}}^{(v)}\alpha_{m\phantom{:1}}^{(v)}}{\alpha^{(v)}} =
\frac{\alpha_{n\phantom{:1}}^{(1)}}{\alpha^{(1)}}\frac{\alpha_{m\phantom{:1}}^{(1)}}{\alpha^{(1)}} -
\frac{\lambda_{n}}{\rho}\left[\frac{\alpha_{n:2}^{(1)}}{\alpha^{(1)}}-\frac{\alpha_{m\phantom{:1}}^{(1)}}{\alpha^{(1)}}\left(\frac{\alpha_{n\phantom{:1}}^{(1)}}{\alpha^{(1)}}\right)^2\right]\,\eta_1^{(v-2)}\nonumber\\[0.0ex]
& \hskip 3.6cm -
\frac{\lambda_{m}}{\rho}\left[\frac{\alpha_{m:2}^{(1)}}{\alpha^{(1)}}-
\frac{\alpha_{n\phantom{:1}}^{(1)}}{\alpha^{(1)}}\left(\frac{\alpha_{m\phantom{:1}}^{(1)}}{\alpha^{(1)}}\right)^2\right]\,\eta_1^{(v-2)}
+ {\rm O}\left(\frac{\lambda^2}{\rho^2}\right)\label{eq:discterm2}\,.
\end{align}
When eqs.~(\ref{eq:discterm1})--(\ref{eq:discterm2}) are put into
eq.~(\ref{eq:gammaik}), an exact cancellation occurs separately among the ${\rm O}(1)$--
and ${\rm O}(\lambda/\rho)$--terms, so we are left with
\begin{equation}
\Gamma^{(v)}_{nm}\ \ \widesim{\ \rho\,\gg\,\max\{\lambda_n,
\lambda_m\}} \ \ \
\frac{\lambda_n}{\rho}\frac{\lambda_m}{\rho}\, {\rm
  O}\left(\frac{\lambda^2}{\rho^2}\right)\,. 
\label{eq:Gammaikexp}
\end{equation} 

\noindent Eqs.~(\ref{eq:Gammaiiexp}) and (\ref{eq:Gammaikexp}) reflect the 
behavior observed in Fig.~\ref{fig:Gammas}.  

\subsection{Asymptotic vanishing of $\eta_k$ from combinatorial arguments}

In order to make the weak truncation expansion effective, we need to
characterize the coefficient functions $\eta_k$ and provide an
algorithmic recipe for their computation. A trivial property,
{\it i.e.} the vanishing of $\eta_k$ as $\rho\to\infty$, can be proved
from purely combinatorial arguments based on
eq.~(\ref{eq:recursion}). This kind of proof nicely follows as a simple application of
the scaling eq.~(\ref{eq:scaling}), yet it gives no clue to the
vanishing rate of $\eta_k$. We begin with the following.
\vskip 0.3cm
\begin{lemma}
Given a set of $m\ge 0$ distinct indices $\{n_1,\dots,n_m\}$ and a corresponding set of
strictly positive multiplicities $\{k_1,\dots,k_m\}$, we have
\begin{equation}
\lim_{\rho\to\infty}\, \rho\,\partial_\rho\alpha^{(v)}_{n_1:k_1\,\ldots\,n_m:k_m}(\rho;\lambda)\, = 0\,.
\label{eq:alphalimit}
\end{equation}
\end{lemma}
\begin{proof}
From eq.~(\ref{eq:recursion}), it follows 
\begin{equation}
\lim_{\rho\to\infty}\,
\rho\,\partial_\rho\alpha^{(v)}_{n_1:k_1\,\ldots\,n_m:k_m} =\,
-\frac{1}{2}\sum_{n=1}^v\lim_{\rho\to\infty}\alpha^{(v)}_{n_1:k_1\,\ldots\,n_m:k_mn}\,
+\,
\frac{1}{2}(v+2k)\lim_{\rho\to\infty}\alpha^{(v)}_{n_1:k_1\,\ldots\,n_m:k_m}\,,
\label{eq:aux1}
\end{equation}
with $k= \sum_{j=1}^m k_j$. The sum index $n$ in eq.~(\ref{eq:aux1})
can either match one of the indices $n_1,\ldots,n_m$ or none of
them. Moreover, we have
\begin{equation}
\lim_{\rho\to\infty}\alpha^{(v)}_{n_1:k_1\,\ldots\,n_m:k_mn} = 
\left\{\begin{array}{ll}
\displaystyle{
\lim_{\rho\to\infty}\alpha^{(v)}_{n_1:k_1\,\ldots\,n_m:k_m}} & \text{if} \quad n\ne n_j\
\ \forall\, j=1,\ldots,m\,; \\[2.0ex]
\displaystyle{
(2k_j+1)\lim_{\rho\to\infty}\alpha^{(v)}_{n_1:k_1\,\ldots\,n_m:k_m}} &
\text{if}\quad n=n_j\ \ \text{for some } j\,;
\end{array}\right.\,,
\end{equation}
as a consequence of the exact factorization of the Gaussian integrals as
$\rho\to\infty$ and the standard formula $\E[z^{2n}] = (2n-1)!!$
($n\ge0$), valid for $z\sim\cN(0,1)$. Therefore, 
\begin{align}
 \lim_{\rho\to\infty}\,
\rho\,\partial_\rho\alpha^{(v)}_{n_1:k_1\,\ldots\,n_m:k_m}\, & =\,
\lim_{\rho\to\infty}\alpha^{(v)}_{n_1:k_1\,\ldots\,n_m:k_m} \cdot\left[
-\frac{v-m}{2}-\frac{1}{2}\sum_{j=1}^m(2k_j+1) + \frac{1}{2}(v+2k)
\right]\nonumber\\[2.0ex]
& = \,\lim_{\rho\to\infty}\alpha^{(v)}_{n_1:k_1\,\ldots\,n_m:k_m}\cdot
\left[\frac{m+2k}{2}- \frac{m+2k}{2}\right] = 0\,.
\end{align}
\end{proof}
\noindent From Lemma 4.1 we easily derive the following.

\vskip 0.4cm

\begin{prop}
As $\rho\to\infty$ all the coefficient functions $\eta_k$ vanish.
%\begin{equation}
%\lim_{\rho\to\infty}\,{\eta^{(v)}_k(\rho;\lambda)} \,=\, 0\,,\qquad k=1,2,\ldots\,.
%\label{eq:zetalimit}
%\end{equation}
\end{prop}
\begin{proof}
Let us define $f_k = (\rho\partial_\rho)^k\alpha^{(v)}$ and $x_k =
\sum_{n_1\ldots n_k=1}^v\alpha^{(v)}_{n_1\dots n_k}$. We first prove by
induction that $f_k$ is a homogeneous linear function of
$x_0,\ldots,x_k$. From eq.~(\ref{eq:recursion}), evaluated at
$k_1=\ldots = k_n=0$, we have indeed
\begin{equation}
\label{eq:reclev0}
 \rho\partial_\rho \alpha^{(v)} = \frac{\,v\,}{2}\,\alpha^{(v)} -
\frac{1}{2}\sum_{k=1}^v\alpha_k^{(v)}\ = \frac{v}{2}x_0 - \frac{1}{2}x_1 = 
f_1\left(x_0,x_1\right)\,.%;\\[1.0ex]
\end{equation}
Now, suppose that $f_{k-1}$ is a homogeneous linear function of
$x_0,\ldots,x_{k-1}$. Then,
\begin{align}
f_k & = (\rho\partial_\rho)^k\alpha^{(v)} =
(\rho\partial_\rho)(\rho\partial_\rho)^{k-1}\alpha^{(v)} =
\rho\partial_\rho f_{k-1}(x_0,\ldots,x_{k-1}) \nonumber\\[2.0ex]
& = f_{k-1}(\rho\partial_\rho x_0,\ldots,\rho\partial_\rho x_{k-1}) = 
f_{k-1}\left(\rho\partial_\rho \alpha^{(v)},\ldots,\sum_{n_1\ldots
  n_{k-1}=1}^v\rho\partial_\rho \alpha^{(v)}_{n_1\ldots n_{k-1}}\right)\,.
\end{align}
The inductive step follows from eq.~(\ref{eq:recursion}) and the
assumed linearity of $f_{k-1}$. Hence, we have
\begin{align}
\lim_{\rho\to\infty}(\rho\partial_\rho)^k\alpha^{(v)} & =
\lim_{\rho\to\infty}f_k(x_0,\ldots,x_k) \nonumber\\[1.0ex]
& \hskip -2.0cm = f_{k-1}\left(\lim_{\rho\to\infty}\rho\partial_\rho\alpha,\ldots,
\sum_{n_1,\ldots,n_{k-1}}\lim_{\rho\to\infty}
\rho\partial_\rho\alpha^{(v)}_{n_1,\ldots,n_k}\right) \, = \, f_{k-1}(0,\ldots,0) \, = \, 0\,,
\end{align}
where the last equality is again a consequence of the homogeneous linearity of
$f_{k-1}$ and the second--to--last one follows from Lemma~4.1. In addition,
we know that (see for instance exercise~13, chap.~6 of ref.~\cite{knuth}) 
\begin{equation}
\rho^k\partial_\rho^{k} = \sum_{j=1}^k(-1)^{k-j}{k\brack
  j}(\rho\partial_\rho)^j\,,
\label{eq:stirlingconn}
\end{equation}
with the symbols ${k\brack j}$ denoting unsigned Stirling numbers of
the first kind. Hence, we conclude
\begin{equation}
\lim_{\rho\to\infty}\eta_k^{(v)} = \lim_{\rho\to\infty}
\frac{1}{\alpha^{(v)}}\sum_{j=1}^k(-1)^{k-j}{k\brack j}(\rho\partial_\rho)^j\alpha^{(v)} =
\sum_{j=1}^k(-1)^{k-j}{k\brack
  j}\lim_{\rho\to\infty}(\rho\partial_\rho)^j\alpha^{(v)} = 0\,.
\end{equation}
\end{proof}

\subsection{Gaussian representation of $\eta_k$}

The above discussion suggests a convenient way to compute the
coefficient functions. We have just seen that $\eta_k$ is a linear
combination of $f_1,\ldots,f_k$. Moreover, $\forall\ j\ge 0$ $f_j$ is a linear
combination of $x_0,\ldots,x_j$. We conclude that $\eta_k$ itself can be represented
as a linear combination of $x_0,\ldots,x_k$. Since we know how to
compute Gaussian integrals with controlled uncertainty, we have a
complete recipe for $\eta_k$, provided we determine the coefficients
of such linear combinations. To  this aim, we concentrate first on the
$f_k$'s. Eq.~(\ref{eq:reclev0}) gives the analytic expression of
$f_1$. By direct calculation we can also derive the expressions
\begin{align}
f_2(x_0,x_1,x_2) & =
\frac{v^2}{4}x_0-\frac{v+1}{2}x_1+\frac{1}{4}x_2\,,\\[0.0ex]
f_3(x_0,x_1,x_2,x_3) & = \frac{v^3}{8}x_0-\frac{3v^2+6v+4}{8}x_1+ \frac{3v+6}{8}x_2-\frac{1}{8}x_3\,,\\[0.0ex]
& \hskip 0.2cm \vdots\nonumber
\end{align}
A generalization is provided by the following
\begin{prop}
For $k\ge1$, we have
\begin{equation}
f_k(x_0,\ldots,x_k) = \sum_{\ell=0}^kd_{k\ell}(v)x_\ell\,,
\label{eq:genfk}
\end{equation}
where the coefficients $d_{k\ell}(v)$ are defined by 
\begin{equation}
d_{k\ell}(v)\, =\, \left\{
\displaystyle{
\begin{array}{cl}
  \displaystyle{ \sum_{t=\ell}^k\frac{(-1)^\ell}{2^t}\phi_{t-\ell}{k\brace
      t}{t\choose \ell}}\,, & k=0,\ldots,\ell\\[4.0ex]
  0\,, & \text{otherwise\,,}
\end{array}}\right.\,,
\label{eq:dkell}
\end{equation}
\begin{equation}
  \phi_k \,=\, \frac{v!!}{(v-2k)!!}\,=\,
    \left\{
    \begin{array}{ll} 1 & k=0 \\[2.0ex]
      (v-2k+2)\,\phi_{k-1} & k\ge 1\end{array}\,\right.\,,
\label{eq:phik}
\end{equation}
and the symbols ${k\brace t}$ denote Stirling numbers of the second kind.
\end{prop}
\begin{proof}
The proof is by induction. We first note that $d_{10}(v) = v/2$ and
$d_{11}(v) = -1/2$. Hence, for $k = 1$ eq.~(\ref{eq:genfk}) agrees
with eq.~(\ref{eq:reclev0}). Now, suppose that $f_k$ is well
represented by eq.~(\ref{eq:genfk}) with $d_{k\ell}(v)$ and $\phi_k$ as in
eqs.~(\ref{eq:dkell}) and (\ref{eq:phik}). Then, from eq.~(\ref{eq:recursion}) it follows
\begin{align}
f_{k+1} & = \, (\rho\partial_\rho)^{k+1}\alpha^{(v)} = \sum_{\ell=0}^k
d_{k\ell}(v)(\rho\partial_\rho)x_\ell = \sum_{\ell=0}^k
d_{k\ell}(v)\sum_{k_1\ldots 
    k_\ell=1}^v\rho\partial_\rho\alpha^{(v)}_{k_1\ldots k_\ell} =
  \nonumber\\[1.0ex]
&  =\, \sum_{\ell=0}^k d_{k\ell}(v)\sum_{k_1\ldots
    k_\ell=1}^v\left\{\frac{1}{2}(v+2\ell)\alpha^{(v)}_{k_1\ldots
    k_\ell}-\frac{1}{2}\sum_{k_{\ell+1}=1}^v\alpha^{(v)}_{k_1\ldots
    k_{\ell+1}}\right\}\nonumber\\[2.0ex]
&  =\, \sum_{\ell=0}^k \left(\frac{v}{2}+\ell\right) d_{k\ell}(v)\sum_{k_1\ldots
    k_\ell=1}^v\alpha^{(v)}_{k_1\ldots k_\ell}
 -\, \frac{1}{2}\sum_{\ell=0}^kd_{k\ell}(v)\sum_{k_1\ldots
    k_{\ell+1}=1}^v\alpha^{(v)}_{k_1\ldots k_{\ell+1}} \nonumber\\[2.0ex]
& =\, \sum_{\ell=0}^k\left[\left(\frac{v}{2}+\ell\right)d_{k\ell}(v) -
    \frac{1}{2}d_{k(\ell-1)}\right]x_\ell\,.
\end{align}
The argument is complete provided we are able to show that
$d_{k\ell}(v)$ fulfills the recurrence
\begin{equation}
  {d_{(k+1)\ell}(v) = \left(\frac{v}{2}+\ell\right)d_{k\ell}(v) -
  \frac{1}{2}d_{k(\ell-1)}(v)}\,.
\end{equation} 
To this aim, it is sufficient to make use of the basic recursive
formulae ${n+1\brace m} = m{n\brace m} + {n\brace m-1}$ and
${n+1\choose m } = {n\choose m} + {n\choose m-1}$. We detail the
algebra for the sake of completeness:
\begin{align}
d_{(k+1)\ell}(v) & = \sum_{t=\ell}^{k+1}\frac{(-1)^\ell}{2^t}\phi_{t-\ell}{k+1 \brace
  t}{t\choose \ell} \nonumber\\[2.0ex]
&  = \sum_{t=\ell}^{k+1}\frac{(-1)^\ell}{2^t}\phi_{t-\ell}\,t{k\brace t}{t\choose \ell} +
\sum_{t=\ell}^{k+1}\frac{(-1)^\ell}{2^t}\phi_{t-\ell}{k\brace t-1}{t\choose \ell}\nonumber\\[2.0ex]
&  = \sum_{t=\ell}^{k+1}\frac{(-1)^\ell}{2^t}\phi_{t-\ell}\,t{k\brace t}{t\choose \ell} + 
\sum_{t=\ell-1}^k\frac{(-1)^\ell}{2^{t+1}}\phi_{t-\ell+1}{k\brace t}{t+1\choose
  \ell}\nonumber
\end{align}
\begin{align}
 & = \sum_{t=\ell}^k\frac{(-1)^\ell}{2^t}\phi_{t-\ell}\,t{k\brace t}{t\choose \ell} +
\sum_{t=\ell-1}^k\frac{(-1)^\ell}{2^{t+1}}\phi_{t-\ell+1}{k\brace t}{t\choose \ell} \nonumber\\[2.0ex]
& + \sum_{t=\ell-1}^k\frac{(-1)^\ell}{2^{t+1}}\phi_{t-\ell+1}{k\brace t}{t\choose \ell-1}
\,=\, \sum_{t=\ell}^k\frac{(-1)^\ell}{2^t}\phi_{t-\ell}\,t{k\brace t}{t\choose
  \ell}\nonumber\\[2.0ex]
& + \sum_{t=\ell-1}^k\frac{(-1)^\ell}{2^{t+1}}[v-2(t-\ell)]\phi_{t-\ell}{k\brace
  t}{t\choose \ell} - \frac{1}{2}d_{k(\ell-1)}(v) = \left(\frac{v}{2}+\ell\right)d_{k\ell}(v) - \frac{1}{2}d_{k(\ell-1)}(v)\,.
\end{align}
\end{proof}
In view of eq.~(\ref{eq:stirlingconn}), it is no surprise that the
coefficients $d_{k\ell}(v)$ embody Stirling 
numbers of the second kind. Recall indeed that Stirling numbers of
the first and the second kind are related to each other by the inversion identity 
\begin{equation}
\sum_{t=0}^{\max\{j,k\}}(-1)^{t-k}{t\brace j}{j\brack t} = \delta_{jk}\,.
\label{eq:invidentity}
\end{equation}
From Proposition~4.2 and eq.~(\ref{eq:invidentity}), it follows
the following proposition.

\begin{prop}
For $k\ge 1$, we have
\begin{equation}
\eta^{(v)}_k(\rho;\lambda) \, = \,
    \sum_{\ell=0}^k\,c_{k\ell}(v)\sum_{k_1\ldots
      k_\ell=1}^v\frac{\alpha^{(v)}_{k_1\ldots
        k_\ell}(\rho;\lambda)}{\alpha^{(v)}(\rho;\lambda)}\,,
\label{eq:etakgauss}
\end{equation}
with the coefficients $c_{k\ell}(v)$ defined as
\begin{equation}
  c_{k\ell}(v) \, = \, \left\{\begin{array}{cl}
  \displaystyle{\frac{(-1)^\ell}{2^k}\frac{v!!}{[v-2(k-\ell)]!!}{k\choose \ell}}\,, & \qquad 0\le \ell\le
  k\,,\\[2.0ex]
  0\,,& \qquad \text{otherwise}\,.
\end{array}\right.
\end{equation}
\end{prop}
\begin{proof}
We have all the necessary ingredients to carry out the proof. Again,
we detail the algebra for the reader's convenience:
\begin{align}
& \eta_k^{(v)} = \frac{1}{\alpha^{(v)}}\rho^k\partial_{\rho}^k\,\alpha^{(v)} =
\frac{1}{\alpha^{(v)}}\sum_{\ell=1}^k(-1)^{k-\ell}{k\brack \ell}(\rho\partial_\rho)^\ell\alpha^{(v)}
%\nonumber\\[2.0ex]
 = \frac{1}{\alpha^{(v)}}\sum_{\ell=1}^k(-1)^{k-\ell}{k\brack \ell}\sum_{m=0}^\ell d_{\ell m}(v)\,x_m\nonumber\\[2.0ex]
& = \frac{1}{\alpha^{(v)}}\sum_{\ell=1}^{k}(-1)^{k-\ell}{k\brack \ell}\sum_{m=0}^\infty\sum_{t=m}^\infty
\frac{(-1)^m}{2^t}\phi_{t-m}{\ell\brace t}{t\choose m}x_m
\nonumber\\[2.0ex]
& = \frac{1}{\alpha^{(v)}}\sum_{m=0}^\infty\sum_{t=m}^\infty
\frac{(-1)^m}{2^t}\phi_{t-m}\sum_{\ell=1}^k(-1)^{k-\ell}{k\brack \ell}{\ell\brace
  t}{t\choose m}x_m \nonumber\\[2.0ex]
& = \frac{1}{\alpha^{(v)}} \sum_{m=0}^\infty \sum_{t=m}^\infty
\frac{(-1)^m}{2^t}\phi_{t-m}\,\delta_{kt}{t\choose m}x_m
 = \frac{1}{\alpha^{(v)}}\sum_{m=0}^\infty
\frac{(-1)^m}{2^k}\phi_{k-m}{k\choose m}x_m \nonumber\\[2.0ex]
& = \frac{1}{\alpha^{(v)}}\sum_{m=0}^k\frac{(-1)^m}{2^k}\phi_{k-m}{k\choose m}x_m\,.
\end{align}
Note that on the second line above we could extend the upper bound of the
sums over $m$ and $t$ from $\ell$ to~$\infty$ thanks to the property ${a\brace b}=0$
if $b>a$. This in turn allowed us to perform the sum exchange on the
third line. By a similar argument, on the last line we could reduce the
upper bound of the sum over $m$ from $\infty$ to~$k$. 
\end{proof}
Obviously, computing eq.~(\ref{eq:etakgauss}) becomes increasingly demanding for
larger values of $k$. Nonetheless, many contributions to the sum on
the \rhs coincide. In particular, all Gaussian integrals with the
same index multiplicities contribute equally, thus we can recast
eq.~(\ref{eq:etakgauss}) to the computationally cheaper expression
\begin{equation}
\eta^{(v)}_k(\rho;\lambda) = \sum_{\ell=0}^kc_{k\ell}(v)\sum_{m_1\ldots m_v=0}^\ell
{\ell\choose
  m_1,\ldots,m_v}\delta_{\ell,m_1+\ldots+m_v}\frac{\alpha^{(v)}_{1:m_1\ldots
    v:m_v}(\rho;\lambda)}{\alpha^{(v)}(\rho;\lambda)}
\,.
\end{equation}

\subsection{Asymptotic sign of $\eta_k$}

In this and next subsection we put the combinatorial approach on hold and work on the
integral representation of the coefficient functions. A first property which turns out
to be essential to the last part of the paper concerns the sign
assumed by $\eta_k$ as $\rho\to\infty$. In regard to this, we state
the following  

\vskip 0.4cm

\begin{prop}
  As $\rho\to\infty$ the sign of $\eta_k$ becomes independent of
  $v$ and $\lambda$. In particular, we have
  \begin{equation}
    \lim_{\rho\to\infty}\sign \eta_k^{(v)}(\rho;\lambda) = (-1)^{k-1}\,.
    \label{eq:asympsign}
  \end{equation}
\end{prop}
\begin{proof}
We first express $\alpha$ in spherical coordinates, {\it i.e.} we perform the
change of integration variable $x = ru$ in eq.~(\ref{eq:alphaints}),
with $r = ||x||$, $u\in\partial\cB_v(1)$ and $\partial\cB_v(1) = \{z\in\RR^v\,:\, ||z||=1\}$
(in the sequel we write $\rd^v x = r^{v-1}\rd r\rd u$; here $\rd u$
embodies the angular part of the spherical Jacobian and the 
differentials of $(v-1)$ angles). Thus, we have
\begin{equation}
\alpha^{(v)}(\rho;\lambda) =
\frac{1}{2\Gamma(v/2)|\Lambda|^{1/2}}\,\M\left[\int_0^{\sqrt{\rho}}\rd
  r\ r^{v-1} \ \exp\left\{-\frac{r^2\cP(u)}{2}\right\}\right]\,,\qquad
\cP(u) = \trans{u}\Lambda^{-1}u\,,
\end{equation}
with $\M$ representing the uniform average operator on
$\partial\cB_v(1)$, namely
\begin{equation}
\M[g] = \frac{\Gamma(v/2)}{2\pi^{v/2}}\int_{\partial\cB_v(1)}\rd
u\ g(u)\,.
\end{equation}
In order to compute $\eta_k$, we differentiate $\alpha$ under the
integral sign. The first derivative evaluates the radial integral at its
upper limit, while the remaining $k-1$ ones distribute according to the chain
rule of differentiation. Explicitly, we have
\begin{align}
\rho^k\partial_\rho^k \alpha^{(v)}(\rho;\lambda) & =
\frac{\rho^k}{2^{v/2}\Gamma(v/2)|\Lambda|^{1/2}}\ \M\left[\partial_\rho^{k-1}\left(\rho^{v/2-1}\exp\left\{-\rho\frac{\cP(u)}{2}\right\}\right)\right]\nonumber
\\[2.0ex]
& = \frac{(-1)^{k-1}\rho^{v/2}}{2^{v/2}\Gamma(v/2)|\Lambda|^{1/2}}\ \M\left[\,\sum_{\ell=0}^{k-1}{k-1
  \choose
  \ell}(-\phi)^{\overline{\ell}}\left(\frac{\rho
    \cP(u)}{2}\right)^{k-1-\ell}\exp\left\{-\frac{\rho \cP(u)}{2}\right\}\right]_{\phi=v/2-1}\nonumber\\[2.0ex]
& =
\frac{(-1)^{k-1}\rho^{v/2}}{2^{v/2}\Gamma(v/2)|\Lambda|^{1/2}}\ \M\left[\cQ_{k-1}\left(\frac{\rho
    \cP(u)}{2},-\phi\right)\exp\left\{-\frac{\rho \cP(u)}{2}\right\}\right]_{\phi=v/2-1}\,.
\label{eq:partderest}
\end{align}
Here $x^{\overline{n}} \equiv x(x+1)\ldots(x+n-1)$ denotes the $n^{\rm
  th}$ raising factorial of $x$ and
\begin{equation}
\cQ_k(x,a) = \sum_{\ell=0}^k{k\choose
  \ell}\,a^{\overline{\ell}}\,x^{k-\ell}\,
\end{equation}
is a polynomial in $x$ of $k^{\rm th}$ degree, differing from a
Newton polynomial for the presence of $a^{\overline{\ell}}$ in
place of $a^\ell$. We note that the coefficient of the leading
term of $Q_k(x,a)$ is ${k\choose 0}a^{\overline{0}}=1$. It follows that
\begin{equation}
\rho^k\partial_\rho^k\alpha^{(v)}(\rho;\lambda)\ \widesim{\ \rho\to\infty\ }\ \frac{(-1)^{k-1}}{\Gamma(\phi+1)|\Lambda|^{1/2}}\left(\frac{\rho}{2}\right)^{\phi+k}\,\M\left[\cP(u)^{k-1}\exp\left\{-\frac{\rho \cP(u)}{2}\right\}\right]
\end{equation}
Eq.~(\ref{eq:asympsign}) follows from the positiveness of $\cP(u)$. 
\end{proof}
\noindent Since $\cP(u)>\lambda_{\rm max}^{-1}$, being $\lambda_{\rm max} =
\max_k\{\lambda_k\}$, we obtain as a by--product an estimate of the
exponential damping of $\eta_k$, namely 
\begin{equation}
\left|\rho^k\partial_\rho^k\alpha^{(v)}(\rho;\lambda)\right| < \frac{\rho^{v/2}}{2^{v/2}\Gamma(v/2)|\Lambda|^{1/2}}\ \left|\M\left[\cQ_{k-1}\left(\frac{\rho
    \cP(u)}{2},-\phi\right)\right]_{\phi=v/2-1}\right|
\re^{-{\rho}/{2\lambda_{\rm max}}}\,.
\label{eq:expdamp}
\end{equation}

\subsection{A convergence estimate for the expansion}

We come back to the issue raised at the beginning of this
section: is the weak truncation expansion convergent? It is not
difficult to see that the answer lies specifically on the behavior of
$\eta_k$ as a function of~$k$. Eq.~(\ref{eq:expdamp}) shows that
the relevant information is brought by $\cQ_{k-1}(\rho\cP(u)/2,-\phi)$, particularly by
its relative minima/maxima. As $k$ increases, the position of the
latter shifts towards larger and larger values of $\rho$, while 
their absolute size increases. In other words, however we choose a
reference scale $\tilde\rho>0$ we always find $\tilde k$ such that
$\argmax\{|\eta_k(\rho;\lambda)|\}>\tilde\rho$ and
$\max\{|\eta_k(\rho;\lambda)|\}>\max\{|\eta_{\tilde
  k}(\rho;\lambda)|\}$ for $k>\tilde 
k$. For this reason, the convergence issue reduces to
quantify the increase rate of $\eta_k$ as a function of $k$.

More quantitatively, we first notice the inequality $\gamma(a,x) <
{x^a}/{a}$. In order to prove this, we observe that a convenient representation of the 
lower incomplete gamma function is provided by  (see for instance sect.~6 of
ref.~\cite{abramowitz}) 
\begin{equation}
\gamma(a,x) \, = \, \frac{1}{a}x^a\re^{-x}\,M(1,1+a,x)\,,
\label{eq:gammaestimate}
\end{equation}
where
\begin{equation}
M(a,b,x) =\,\!_1F_1(a;b;x) \, = \, \sum_{n=0}^\infty\frac{a^{\overline
    n}}{b^{\overline{n}}}\frac{x^n}{n!}
\end{equation}
is the confluent hypergeometric function originally introduced by Kummer. Since
$1^{\overline{n}}=n!$ and $(1+a)^{\overline{n}} =
(1+a)(2+a)\ldots(n+a)>n!$ if $a>0$, it follows 
\begin{equation}
M(1,1+a,x) =
\sum_{n=0}^\infty\frac{1^{\overline{n}}}{(1+a)^{\overline{n}}}\frac{x^n}{n!} =
\sum_{n=0}^\infty\frac{x^n}{(1+a)^{\overline{n}}}<\sum_{n=0}^\infty
\frac{x^n}{n!} = \re^x\,.
\end{equation}
We can use eq.~(\ref{eq:gammafunc}) and the above inequality to
establish an upper bound to the 1--dimensional Gaussian integrals,
namely 
\begin{equation}
\alpha^{(1)}_{n:k}(\rho;\lambda_k) < \frac{1}{\sqrt{2\pi}k}\left(\frac{\rho}{\lambda_n}\right)^{k+1/2}\,.
\label{eq:oneindgammaest}
\end{equation}
Now, let us denote by $\cX_{n:k}$ the weak truncation expansion of
$\alpha_{n:k}$, {\it i.e.}
\begin{equation}
\cX_{n:k}^{(v)}(\rho;\lambda) = \alpha^{(v-1)}(\rho;\lambda_{(n)})
\sum_{p=0}^\infty\frac{(-1)^p}{p!}\left(\frac{\lambda_n}{\rho}\right)^p \alpha_{n:(k+p)}^{(1)}(\rho;\lambda_n)\eta_p^{(v-1)}(\rho;\lambda_{(n)})\,.
\end{equation}
In view of eq.~(\ref{eq:oneindgammaest}), an absolute estimate to $\cX_{n:k}$ is given by
\begin{align}
|\cX^{(v)}_{n:k}(\rho;\lambda)| & \,<\, \sum_{p=0}^\infty
\frac{1}{p!}\left(\frac{\lambda_n}{\rho}\right)^p\alpha^{(1)}_{n:(k+p)}(\rho;\lambda_n)\alpha^{(v-1)}(\rho;\lambda_{(n)})\left|\eta_p^{(v-1)}(\rho;\lambda_{(n)})\right|\nonumber\\[2.0ex]
& \,<\,
\frac{1}{\sqrt{2\pi}}\left(\frac{\rho}{\lambda_n}\right)^{k+1/2}\sum_{p=0}^\infty
\frac{1}{p!}\frac{1}{(p+k)}\left|\rho^p\partial_\rho^p\alpha^{(v-1)}(\rho;\lambda_{(k)})\right|\,.
\label{eq:Xnfirstest}
\end{align}

\noindent From the above inequality we see that a less than factorial growth of
$\eta_p$ with $p$ would make the \rhs of eq.~(\ref{eq:Xnfirstest})
convergent. In order to 
estimate $|\rho^p\partial_\rho^p\alpha|$, we make use of
eq.~(\ref{eq:partderest}). Since $\partial\cB_v(1)$ is a compact
domain, we can get rid of the angular average by defining
\begin{equation}
u^* =
\argmax_{u\,\in\,\partial\cB_v(1)}\left\{\left|\cQ_{p-1}\left(\frac{\rho
\cP(u)}{2},-\phi\right)\right|\exp\left\{-\frac{\rho
    \cP(u)}{2}\right\}\right\}\,,
\end{equation}
whence it follows
\begin{align}
\left|\rho^p\partial_\rho^p\alpha^{(v)}(\rho;\lambda)\right| & \le
\frac{\rho^{v/2}}{2^{v/2}\Gamma(v/2)|\Lambda|^{1/2}}\,\left|\cQ_{p-1}\left(\frac{\rho 
  \cP(u^*)}{2},-\phi\right)\right|\exp\left\{-\frac{\rho
    \cP(u^*)}{2}\right\}\,.
\label{eq:secderest}
\end{align}

\noindent We have already noticed that $\cP(u)>\lambda_{\rm max}^{-1}$ for all $u\in\partial\cB_v(1)$. 
If we also consider that $|\Lambda|^{1/2}>\lambda_{\rm min}^{v/2}$,
being $\lambda_{\rm min} = \min_k\{\lambda_k\}$, then we have $[\cP(u)^{v/2}|\Lambda|^{1/2}]^{-1}<(\lambda_{\rm
  max}/\lambda_{\rm min})^{1/2}$. Multiplying and dividing the
\rhs of eq.~(\ref{eq:secderest}) by $\cP(u^*)^{v/2}$ leads to
\begin{align}
\left|\rho^p\partial_\rho^p\alpha^{(v)}(\rho;\lambda)\right|&
\,<\,\frac{1}{\Gamma(v/2)}\left(\frac{\lambda_{\rm max}}{\lambda_{\rm min}}\right)^{v/2}\left(\frac{\rho
\cP(u^*)}{2}\right)^{v/2}\left|\cQ_{p-1}\left(\frac{\rho
  \cP(u^*)}{2},-\phi\right)\right|\exp\left\{-\frac{\rho
  \cP(u*)}{2}\right\}\nonumber\\[2.0ex]
& \, < \, \frac{1}{\Gamma(v/2)}\left(\frac{\lambda_{\rm
    max}}{\lambda_{\rm min}}\right)^{v/2}\max_{x\,\in\,\RR_+}\left\{x^{v/2}\left|\cQ_{p-1}(x,-\phi)\right|\re^{-x}\right\}\,.
\end{align}
Accordingly, we obtain the estimate
%% \begin{align}
%% |\cX_{n:k}^{(v)}(\rho;\lambda)| & \,<\,\frac{1}{\sqrt{2\pi}\Gamma((v-1)/2)}\left(\frac{\rho}{\lambda_n}\right)^{k+1/2}\left(\frac{\tilde\lambda_{\rm
%%     max}}{\tilde\lambda_{\rm min}}\right)^{(v-1)/2}\nonumber\\[2.0ex]
%% & \hskip 2.0cm \cdot \sum_{p=1}^\infty\frac{1}{p!}\frac{1}{p}\max_{x\,\in\,\RR_+}\left\{x^{(v-1)/2}\left|\cQ_{p-1}(x,-\phi^*)\right|\re^{-x}\right\}\,,
%% \end{align}
\begin{align}
|\cX_{n:k}^{(v)}(\rho;\lambda)| & \,<\,\frac{1}{\sqrt{2\pi}\Gamma((v-1)/2)}\left(\frac{\tilde\lambda_{\rm
    max}}{\tilde\lambda_{\rm min}}\right)^{(v-1)/2}\left(\frac{\rho}{\lambda_n}\right)^{k+1/2} \sum_{p=1}^\infty\frac{\cC^{(v)}(p)}{p}\,,
\label{eq:finest}
\end{align}
where we have set $\tilde\lambda_{\rm min} = \min_{j\ne n}\{\lambda_j\}$, $\tilde\lambda_{\rm max} =
\max_{j\ne n}\{\lambda_{j}\}$ and
%\begin{align}
%\cC^{(v)}(p) & =
%\frac{1}{p!}\max_{x\,\in\,\RR_+}\left\{x^{(v-1)/2}\left|\cQ_{p-1}(x,-\phi^*)\right|\re^{-x}\right\}\nonumber\\[2.0ex]
%& = \max_{x\,\in\,\RR_+}\left\{\left|\sum_{\ell=0}^{p-1}\frac{\,(-\phi^*)^{\overline\ell}\,x^{p-\ell-\phi^*}}{\ell!(p-1-\ell)!}\right|\,\re^{-x}\right\}_{\phi^*= (v-3%)/2}\,.
%\end{align}
\begin{equation}
\cC^{(v)}(p) = \frac{1}{p}\max_{x\,\in\,\RR_+}\left\{\left|\sum_{\ell=0}^{p-1}\frac{\,(-\phi^*)^{\overline\ell}\,x^{p-\ell-\phi^*}}{\ell!(p-1-\ell)!}\right|\,\re^{-x}\right\}_{\phi^*= (v-3)/2}\,.
\end{equation}

\noindent It will be observed that in order to arrive at eq.~(\ref{eq:finest}),
we have gone through a rather long inequality 
chain, so it is not clear whether the resulting upper bound is
finite. For the sum on the \rhs to be convergent, it suffices that
$\exists\,\epsilon>0$, $A>0$ and $p_0$ such that
$\cC(p)<Ap^{-\epsilon}\equiv m(p;A,\epsilon)$ for $p>p_0$. If this holds true, then we have
$\sum_{p=1}^{\infty}\cC(p)/p< A\zeta(1+\epsilon) + \text{const.}$, 
where $\zeta$ denotes the Riemann zeta function. In order to
evaluate $\cC(p)$ analytically, we need to solve a polynomial equation
of degree $(p-\phi^*)$ for $p\gg 1$. An alternative approach is to compute
$\cC(p)$ numerically for a set of sufficiently large values of $p$ and
then try to fit data to a model such as $m(p;A,\epsilon)$, with
the parameters $A$ and $\epsilon$ depending in general on $v$. In
Fig.~\ref{fig:Cp} (right) we plot numerical determinations of $\cC(p)$
for $v=2,\ldots,6$. We observe that 
$\cC(p)$ is monotonic decreasing for $v\le 5$ and monotonic increasing
for $v\ge 6$. In Fig.~\ref{fig:Cp} (left) we report the values of the
fitted parameters $A$ and $\epsilon$ together with the corresponding
$\chi^2$--values (the range chosen for all fits is
$p\in[50,100]$). Our numerical experiments suggest that the weak 
truncation expansion converges uniformly in $\rho$ at least for $v\le
5$. Nevertheless, the argument is not conclusive, since it assumes
that we can extrapolate the fitting model to $p\to\infty$, which is
not mathematically rigorous...  

\begin{center}
  \begin{figure}
  \vskip 1.2cm
    \centering
    \hskip -5.0cm
    \subfigure{ 
      \hskip 5.0cm\begin{tabular}{c|ccc}
        %\\[-8.0cm]
        \hline\hline\\[-2ex]
        $v$ & $A$ & $\phantom{-}\epsilon$ & $\chi^2/\text{ndf}$ \\
        \hline \\[-1.5ex]
        2 & 0.522 & $\phantom{-}0.734$ & $1.08\times 10^{-11}$ \\
        3 & 0.396 & $\phantom{-}0.499$ & $9.86\times 10^{-13}$ \\
        4 & 0.361 & $\phantom{-}0.278$ & $5.63\times 10^{-10}$ \\
        5 & 0.375 & $\phantom{-}0.068$ & $2.16\times 10^{-8\phantom{1}}$ \\
        6 & 0.227 & $-0.309$ & $3.30\times 10^{-7\phantom{1}}$ \\[0.3ex] 
        \hline\hline
      \end{tabular}
    }
    \subfigure{\hskip 7.0cm}
    \vskip -4.2cm
    \hskip 2.0cm\subfigure{\hskip 4.0cm\includegraphics[width=0.40\textwidth]{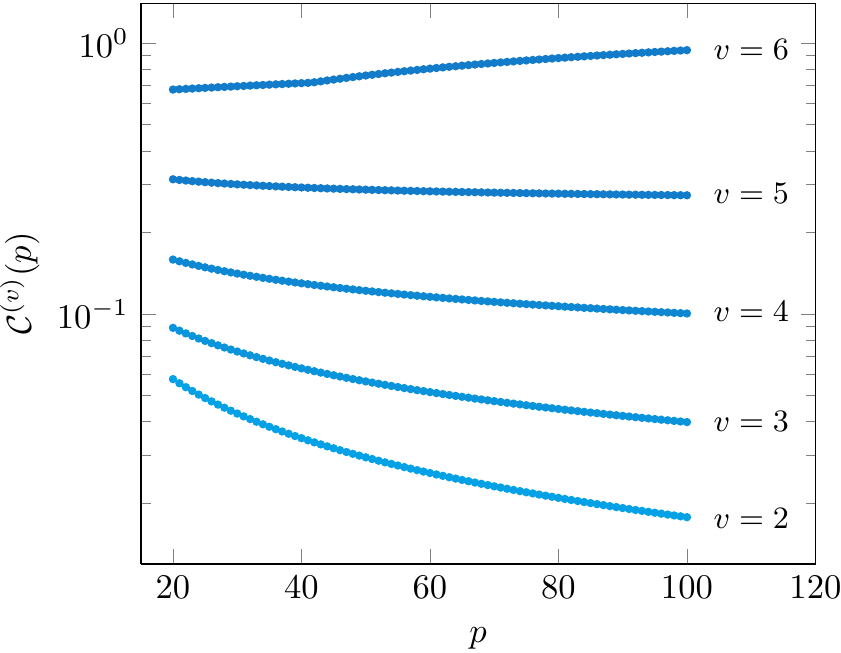}}
    \caption{Numerical computation of $\cC^{(v)}(p)$ for
      $v=2,\ldots,6$. We fit each curve on the right to a model $m(p;A,\epsilon) =
      Ap^{-\epsilon}$, with $A$ and $\epsilon$ depending in general on
      $v$. Fitted parameters are reported on the left,
      together with the corresponding $\chi^2$--values.\label{fig:Cp}}
  \end{figure}
  \vskip -1.0cm
\end{center}

\section{Variance reduction in the regime of weak truncation}

If we look at eq.~(\ref{eq:weakvar}), we see that the next--to--leading
contribution to the rescaled variance $\Gamma_{nn}$ is the sum of
a few ratios of 1--dimensional Gaussian integrals, all proportional to
$\eta_1$. The subsequent terms of the expansion have an increasingly
complex structure. Each power of $\lambda_n/\rho$ couples to several
products of coefficient functions $\eta_k$, always combined so as to
 give the correct power counting. If $f(\rho;\lambda)$ is a
generic observable, its weak truncation expansion reads
\begin{equation}
f(\rho;\lambda) = \sum_{q=0}^{\infty}
(-1)^q\left(\frac{\lambda_n}{\rho}\right)^q\sum_{\ue\in\cS_q^q}\,
\Xi_f^{(q;\ue)}(\rho,\lambda)\
\eta_0^{(v-1)}(\rho;\lambda_{(n)})^{e_0}\ldots\eta_q^{(v-1)}(\rho;\lambda_{(n)})^{e_q}\,.
\label{eq:genwte}
\end{equation}
We denote by $\ue = \{e_0,\ldots,e_q\}$ the exponents
of  $\eta_0,\ldots,\eta_q$ and by
$\cS_q^m$ the set of all possible $\ue$'s corresponding to an overall power
counting $m$, namely
\begin{equation}
\cS_q^m = \left\{\,\ue\in\NN_0^{q+1}:\ \ \cP_q(\ue)=m\,\right\}\,,
\end{equation}
with the power counting function $\cP_q(\ue)$  defined as
\begin{equation}
\cP_q(\ue) = \sum_{k=1}^{q} k\, e_k\,.
\end{equation}
Note that if $m<q$ and $\ue\in\cS_q^m$, then
$e_{m+1}=\ldots=e_q=0$. In this case, we interpret $\Xi^{(m;\ue)}_f$
as $\Xi_f^{(m;\{e_0,\ldots,e_m\})}$. Recall also that since $\eta_0=1$, 
$e_0$ never contributes to the power counting. For later convenience, we 
define $\ove=\{e_1,\ldots,e_q\}$. If $\ue\in\cS_q^m$, with abuse of
notation we also write $\ove\in\cS_q^m$. Strictly speaking, the
presence of $\eta_0$ in eq.~(\ref{eq:genwte}) is necessary in order to
properly take into account the leading order of the expansion, namely
\begin{equation}
\lim_{\rho\to\infty}f(\rho;\lambda) = \,\lim_{\rho\to\infty}\sum_{e_0=0}^{\infty}\Xi_f^{(0;\{e_0\})}(\rho;\lambda)\,.
\end{equation}
The sum over $e_0$ in eq.~(\ref{eq:genwte}) extends in principle from 0 to
$\infty$. We use $e_0$ to enumerate all contributions proportional to
$\eta_1^{e_1}\ldots\eta_q^{e_q}$. Accordingly, the information concerning the
maximum value taken by $e_0$ is hidden within $\Xi_f^{(q;\ue)}$. We finally stress that  
$\Xi_f^{(q;\ue)}$ depends in general upon all the components of
$\lambda$, yet in the specific cases of $\Gamma_{nn}$ and $\Delta_n$
(which are the ones we are interested in) it depends only upon $\lambda_n$.  

Suppose now that $f(\rho;\lambda)$ and $g(\rho;\lambda)$ are two observables, which we
expand according to eq.~(\ref{eq:genwte}). It is not difficult to
prove that the convolution rules needed to obtain the expansion
of the algebraic combinations $f+g$ and  $f\cdot g$ are similar to
the Fourier ones. Specifically, we have 
\begin{equation}
\Xi_{f+g}^{(q;\ue)} = \Xi_f^{(q;\ue)}  + \Xi_g^{(q;\ue)}\,,
\label{eq:xisum}
\end{equation}
and
\begin{equation}
\Xi_{f\cdot g}^{(q;\ue)} = \sum_{\ell,m=0}^{q}\,
\delta_{q,\ell+m} \,\sum_{\uc\in\cS_q^\ell}\,\sum_{\ud\in\cS_q^m}\, \delta_{\ue,\uc+\ud} \ \Xi_f^{(\ell;\uc)}\Xi_g^{(m;\ud)}\,,
\label{eq:xiprod}
\end{equation}
where $\delta_{\uc,\ud} = \prod_{k=0}^{q}\delta_{c_k,d_k}$ is a
vector generalization of
the Kronecker symbol. Eqs.~(\ref{eq:xisum}) and (\ref{eq:xiprod}) are
sufficient to derive the weak truncation expansion of $\Delta_n$ to
all orders and to consequently prove our main result, represented by
the following 

\begin{center}
  \fcolorbox{blue}{shadecolor}
            {
              \parbox{0.85\textwidth}
                     {
                       \begin{theo}
                         As $\rho\to\infty$, the sign of the coefficient function
                         $\Xi_{\Delta_n}^{(q;\ue)}(\rho;\lambda_n)$ is given by
                         \begin{equation}
                           \lim_{\rho\to\infty}\sum_{\,e_0}\Xi^{(q;\ue)}_{\Delta_n}(\rho;\lambda_n) = (-1)^{\sum_{k=1}^qe_k-1}\,.
                         \end{equation}
                         This, in conjunction with eq.~(\ref{eq:asympsign}), implies that all
                         terms of the weak truncation expansion of $\Delta_n$ become negative
                         at sufficiently large $\rho$. 
                       \end{theo}
                     }
            }
\end{center}

\begin{proof}
We proceed in subsequent steps. First of all, we observe that
$\Delta_n$ can be written as
\begin{equation}
\Delta_n(\rho;\lambda) =
\frac{\lambda_n^2}{\rho^2}\,\fDn\cdot\fDd\,;\qquad 
\left\{\begin{array}{l}
\fDn = \alpha^{(v)}_{nn}\,\alpha^{(v)} - \alpha_n^{(v)}\alpha_n^{(v)} - 2\,\alpha_n^{(v)}\alpha^{(v)}\,,\\[2.0ex]
\fDd = [\alpha^{(v)}]^{-2}\,.
\end{array}\right.
\end{equation}
To work out $\fDn$, we first review the expansion of $\alpha_{n:k}$. Explicitly, we have
\begin{align}
\alpha^{(v)}_{n:k} & = \sum_{q=0}^\infty
\frac{(-1)^q}{q!}\left(\frac{\lambda_n}{\rho}\right)^q\alpha^{(1)}_{n:(k+q)}\alpha^{(v-1)}\eta_q^{(v-1)}\nonumber\\[2.0ex]
& =
\sum_{q=0}^\infty\frac{(-1)^q}{q!}\left(\frac{\lambda_n}{\rho}\right)^q\alpha^{(1)}_{n:(k+q)}\alpha^{(v-1)}\sum_{\ue\in\cS_q^q}\,\left(\prod_{i=0}^{q-1}\delta_{e_i,0}\right)\,\delta_{e_q,1}\,[\eta^{(v-1)}_0]^{e_0}\ldots[\eta^{(v-1)}_q]^{e_q}\,.
\end{align}
Hence, we deduce
\begin{equation}
\Xi^{(q;\ue)}_{\alpha_{n:k}} = \frac{1}{q!}\alpha^{(1)}_{n:(k+q)}\alpha^{(v-1)}\,\left(\prod_{i=0}^{q-1}\delta_{e_i,0}\right)\,\delta_{e_q,1}\,.
\end{equation}
This allows us to derive the expansion of the product
$\alpha_{n:r}\alpha_{n:s}$. From eq.~(\ref{eq:xiprod}) it follows
\begin{align}
\Xi^{(q;\ue)}_{\alpha_{n:r}\alpha_{n:s}} & =
\sum_{\ell,m=0}^{q}\frac{\delta_{q,\ell+m}}{\ell!\,m!}\,\alpha^{(1)}_{n:(r+\ell)}\alpha^{(1)}_{n:(s+m)}[\alpha^{(v-1)}]^2\nonumber\\[2.0ex]
& \cdot\, \sum_{\uc\in\cS_q^\ell}\,\sum_{\ud\in\cS_q^m}\delta_{\ue,\uc+\ud}\ \left(\prod_{i=0}^{\ell-1}\delta_{c_i,0}\right)\left(\prod_{k=0}^{m-1}\delta_{d_k,0}\right)\delta_{c_\ell,1}\delta_{d_m,1}\,.
\end{align}
The inner sums can be performed exactly. Indeed, thanks to the
Kronecker symbols $\delta_{c_\ell,1}$ and $\delta_{d_m,1}$, 
non--vanishing contributions group according to whether $\ell=m$ or
$\ell\ne m$, namely
\begin{align}
\Xi^{(q;\ue)}_{\alpha_{n:r}\cdot\alpha_{n:s}} =
[\alpha^{(v-1)}]^2\, & \sum_{\ell,m=0}^{q}
  \frac{\delta_{q,\ell+m}}{\ell!\,m!}\,\alpha_{n:(r+\ell)}^{(1)}\alpha_{n:(s+m)}^{(1)}\nonumber\\[2.0ex]
& \cdot \left[\delta_{\ell,m}\delta_{e_\ell,2}\prod_{i\ne
    \ell}^{0\ldots q}\delta_{e_i,0}+(1-\delta_{\ell,m})\delta_{e_\ell,1}\delta_{e_m,1}\prod_{i\ne
  \ell,m}^{0\ldots q}\delta_{e_i,0}\right]\,.
\label{eq:wteprod}
\end{align}
We see that the Kronecker symbol $\delta_{q,\ell+m}$ makes both
groups of terms vanish unless $\cP_q(\ue) = q$, as intuitively
understood. Conversely, the only elements $\ue\in\cS_q^q$ which result in
a non--vanishing coefficient function
$\Xi^{(q,\ue)}_{\alpha_{n:r}\alpha_{n:s}}$ are either those where two
different exponents $e_\ell$, $e_m$ equal one (with $\ell+m=q$) while the others vanish, or
those where $e_{q/2}=2$ and $e_i=0$ for $i\ne q/2$ (of
course the latter contribute only when $q$ is even). From
eq.~(\ref{eq:wteprod}), we immediately obtain
\begin{align}
\Xi^{(q;\ue)}_{\fDn}  =
   [\alpha^{(v-1)}]^2 & \,\sum_{\ell,m=0}^q\frac{\delta_{q,\ell+m}}{\ell!\,m!}\left(\alpha^{(1)}_{n:(\ell+2)}\,\alpha^{(1)}_{n:m}
     - \alpha^{(1)}_{n:(\ell+1)}\alpha^{(1)}_{n:(m+1)} -
     2\alpha^{(1)}_{n:(\ell+1)}\alpha^{(1)}_{n:m}\,\right)\nonumber\\[2.0ex]
& \cdot  \left[\delta_{\ell,m}\delta_{e_\ell,2}\prod_{i\ne
    \ell}^{0\ldots q}\delta_{e_i,0}+(1-\delta_{\ell,m})\delta_{e_\ell,1}\delta_{e_m,1}\prod_{i\ne
  \ell,m}^{0\ldots q}\delta_{e_i,0}\right]\,.
\end{align}
The above expression depends upon $\rho$ essentially via the integrals in
parentheses (the overall factor $[\alpha^{(v-1)}]^2$ is irrelevant to
our aims). Since $\alpha_{n:r}\to (2r-1)!!$ as 
$\rho\to\infty$, we have $(\alpha_{n:(\ell+2)}\,\alpha_{n:m} 
     - \alpha_{n:(\ell+1)}\alpha_{n:(m+1)} -
     2\alpha_{n:(\ell+1)}\alpha_{n:m})\to
     (\ell-m)(2\ell+1)(2\ell-1)!!(2m-1)!!$. In particular, this
     quantity vanishes for $\ell=m$, thus making the first term
     in square brackets never contribute as $\rho\to\infty$. A little
     additional algebra yields
\begin{equation}
\lim_{\rho\to\infty}\Xi^{(q;\ue)}_{\fDn}(\rho;\lambda) =
4\sum_{\ell=0}^q\sum_{m=0}^{\ell-1}\delta_{q,\ell+m}(\ell-m)^2\frac{(2\ell-1)!!}{\ell!}\frac{(2m-1)!!}{m!}\delta_{e_\ell,1}\delta_{e_m,1}\prod_{i\ne
\ell,m}^{0\ldots q}\delta_{e_i,0}\,.
\label{eq:xifDn}
\end{equation}
We notice that the \rhs of eq.~(\ref{eq:xifDn}) vanishes always for $e_0\ge2$, but
not necessarily for $e_0=0$ or $e_0=1$. 

As a second step, we work out $\fDd$. To this aim, we first need to evaluate
$\Xi^{(q;\ue)}_{\alpha^{-1}}$. As already
done in sect.~4, we make use of the Taylor series $(1+x)^{-1} = \sum_{p=0}^\infty
(-1)^px^p$. From
\begin{equation}
\alpha^{(v)} = \alpha^{(1)}\alpha^{(v-1)}\cdot\left[1+\sum_{q=1}^\infty\frac{(-1)^q}{q!}\left(\frac{\lambda_n}{\rho}\right)^q\frac{\alpha^{(1)}_{n:q}}{\alpha^{(1)}}\,\eta_q^{(v-1)}\right]\,,
\end{equation}
it follows
\begin{align}
[\alpha^{(v)}]^{-1} & =
[\alpha^{(1)}\alpha^{(v-1)}]^{-1}\sum_{p=0}^{\infty}(-1)^p\sum_{\ell_1\ldots
\ell_p=1}^\infty \frac{(-1)^{\ell_1+\ldots+\ell_p}}{\ell_1!\ldots
  \ell_p!}\left(\frac{\lambda_n}{\rho}\right)^{\ell_1+\ldots+\ell_p}\frac{\alpha^{(1)}_{n:\ell_1}\ldots
\alpha^{(1)}_{n:\ell_p}}{[\alpha^{1}]^p}\eta_{\ell_1}^{(v-1)}\ldots\eta_{\ell_p}^{(v-1)}\nonumber\\[2.0ex]
& = [\alpha^{(1)}\alpha^{(v-1)}]^{-1}\sum_{q=0}^\infty
(-1)^q\left(\frac{\lambda_n}{\rho}\right)^q
\sum_{p=0}^q(-1)^p\sum_{\ell_1\ldots \ell_p=1}^\infty \delta_{q,\ell_1+\ldots+\ell_p}\prod_{j=1}^p\left[\frac{1}{\ell_j!}\frac{\alpha^{(1)}_{n:\ell_j}}{\alpha^{(1)}}\eta_{\ell_j}^{(v-1)}\right]\,.
\end{align}
On the second line we have reduced the upper limit of the sum over $p$ from
$\infty$ to $q$. The reason is that
$\delta_{q,\ell_1+\ldots+\ell_p}=0$ for $q<p$, owing to
$\ell_1+\ldots+\ell_p\ge p$. On expanding the sums over 
$\ell_1,\ldots,\ell_p$, we see that all terms proportional to
$\eta_1^{e_1}\ldots\eta_q^{e_q}$ for some $\ove$
coincide. Since permutations of $\ell_1,\ldots,\ell_p$
corresponding to the same $\ove$ give all the same contribution, the latter
turns out to be multiplied by an overall numerical factor which is at most
$p!$. Of course, permutations of equal indices
contribute only once. Therefore, a correct counting of that factor
amounts to the multinomial coefficient $p!/(e_1!\ldots e_q!)$, with the
constraint $\sum e_j = p$. In other words, we have 
\begin{equation}
[\alpha^{(v)}]^{-1} =
[\alpha^{(1)}\alpha^{(v-1)}]^{-1}\sum_{q=0}^\infty
(-1)^q\left(\frac{\lambda_n}{\rho}\right)^q\sum_{\ue\in\cS_q^q}(-1)^{\sum
e_k}\delta_{e_0,0}{\sum e_k\choose e_0,\ldots,e_k}\prod_{j=0}^q\left[\frac{1}{j!}\frac{\alpha^{(1)}_{n:j}}{\alpha^{(1)}}\eta_j\right]^{e_j}\,,
\end{equation}
and consequently
\begin{equation}
\Xi^{(q;\ue)}_{\alpha^{-1}} = (-1)^{\sum_{k=1}^q
e_k}\,\delta_{e_0,0}{\sum e_k\choose
  e_0,\ldots,e_k}\prod_{j=0}^q\left[\frac{1}{j!}\frac{\alpha^{(1)}_{n:j}}{\alpha^{(1)}}\right]^{e_j}
[\alpha^{(1)}\alpha^{(v-1)}]^{-1}\,.
\label{eq:xiam1}
\end{equation}
\vskip -0.3cm
\noindent Now, we obtain $\Xi_{\fDd}^{(q;\ue)}$ by convolving eq.~(\ref{eq:xiam1}) with itself. This yields
\begin{equation}
\Xi^{(q;\ue)}_{\fDd} = \left\{(-1)^{\sum_{k=1}^q
  e_k}\,\delta_{e_0,0}\prod_{j=1}^q\left[\frac{1}{j!}\frac{\alpha^{(1)}_{n:j}}{\alpha^{(1)}}\right]^{e_j}\Psi^{(q;\ove)}\right\}[\alpha^{(1)}\alpha^{(v-1)}]^{-2}\,,
\label{eq:xifDdnolim}
\end{equation}
with the coefficient $\Psi^{(p;\ove)}$ defined by
\begin{equation}
\Psi^{(p;\ove)} = \sum_{\ell,m=0}^p\delta_{p,\ell+m}\sum_{\ovc\in\cS_q^\ell}\sum_{\ovd\in\cS_q^m}\delta_{\ove,\ovc+\ovd}\,{\sum_{k=1}^q
  c_k\choose c_1,\ldots,c_q}{\sum_{k=1}^q d_k\choose
  d_1,\ldots,d_q}\,,\qquad |\ove|=q\,.
\end{equation}
From eq.~(\ref{eq:xifDdnolim}), it follows
\begin{equation}
\lim_{\rho\to\infty}\Xi^{(q;\ue)}_{\fDd} = (-1)^{\sum_{k=1}^q e_k}\,\delta_{e_0,0}\prod_{j=1}^q\left[\frac{(2j-1)!!}{j!}\right]^{e_j}\Psi^{(q;\ove)}\,.
\label{eq:xifDd}
\end{equation}
\vskip -0.3cm
As a third step, we convolve eqs.~(\ref{eq:xifDn}) and
(\ref{eq:xifDd}). In this way we obtain $\Xi^{(q;\ue)}_{\Delta_n}$
directly in the limit $\rho\to\infty$. The algebra is just a little
bit intricate, so we detail it. First of all,
\begin{align}
\lim_{\rho\to\infty}\Xi^{(q;\ue)}_{\Delta_n} & =
4\sum_{\ell,m=0}^q\delta_{q,\ell+m}\sum_{\uc\in\cS_q^\ell}\sum_{\ud\in\cS_q^m}\delta_{e_0,c_0}\delta_{\ove,\ovc+\ovd}\left[\sum_{r=0}^\ell\sum_{s=s}^{r-1}\delta_{n,r+s}(r-s)^2\frac{(2r-1)!!}{r!}\frac{(2s-1)!!}{s!}\right.\nonumber\\[2.0ex]
& \cdot \left.\delta_{c_r,1}\delta_{c_s,1}\prod_{i\ne
r,s}^{0..\ell}\delta_{c_i,0}\right]\cdot\left\{(-1)^{\sum_{k=1}^m
  d_k}\prod_{k=1}^m\left[\frac{(2k-1)!!}{k!}\right]^{d_k}\Psi^{(m;\ovd)}\right\}\,.
\end{align}
Owing to the Kronecker symbols, we pay no price if we  introduce a
factor of $(-1)^{\sum_{k=0}^q c_k}\equiv 1$ within square brackets. For the same
reason, we can also insert additional factors of
$[(2k-1)!!/k!]^{c_k}\equiv 1$ for $k\ne r,s$ without pay. Hence,
\begin{align}
\lim_{\rho\to\infty}\Xi^{(q;\ue)}_{\Delta_n} & = 4 (-1)^{\sum_{k=1}^q
  e_k}
\prod_{k=1}^q\left[\frac{(2k-1)!!}{k!}\right]^{e_k}\sum_{\ell,m=0}^q\delta_{q,\ell+m}\sum_{\uc\in\cS_q^\ell}\sum_{\ovd\in\cS_q^m}(-1)^{c_0}\delta_{e_0,c_0}\delta_{\ove,\ovc+\ovd}\nonumber\\[0.0ex]
&\cdot \sum_{r=0}^\ell\sum_{s=0}^{r-1}\delta_{\ell,r+s}(r-s)^2
\delta_{c_r,1}\delta_{c_s,1}\prod_{i\ne
  r,s}^{0\ldots\ell}\delta_{c_i,0}\,\Psi^{(m;\ovd)}\nonumber\\[0.0ex]
& = 4(-1)^{\sum_{k=1}^q
  e_k}\delta_{\cP_q(\ue),q}\prod_{k=1}^q\left[\frac{(2k-1)!!}{k!}\right]^{c_k}\sum_{\ell,m=0}^q\delta_{q,\ell+m}\sum_{\uc\in\cS_q^\ell}
(-1)^{c_0}\delta_{e_0,c_0}\nonumber\\[0.0ex]
&
\cdot\sum_{r=0}^\ell\sum_{s=0}^{r-1}\delta_{\ell,r+s}(r-s)^2\delta_{c_r,1}\delta_{c_s,1}\prod_{i\ne
r,s}^{0\ldots\ell}\delta_{c_i,0}\,\Psi^{(m;\ove-\ovc)}\,\theta_{\ue,\uc}\,,
\label{eq:finalDeltan}
\end{align}
with $\theta_{\ove,\ovc}=\prod_{i=1}^{q}\theta_{e_i,c_i}$ being a
vector generalization of the Heaviside function
\begin{equation}
\theta_{a,b} = \left\{\begin{array}{ll} 1 & \text{if } a\ge b\,,\\[2.0ex] 0
& \text{otherwise}\,.\end{array}\right.
\end{equation}
From eq.~(\ref{eq:finalDeltan}) it follows
\begin{equation}
\lim_{\rho\to\infty}\sum_{e_0=0}^\infty \Xi^{(q;\ue)}_{\Delta_n} = 4(-1)^{\sum_{k=1}^qe_k}\delta_{\cP_q(\ue),q}\prod_{k=1}^q\left[\frac{(2k-1)!!}{k!}\right]^{e_k}\left\{\Omega_0^{(q;\ove)}-\Omega_1^{(q;\ove)}\right\}\,,
\end{equation}
with
\begin{equation}
\Omega_0^{(q;\ove)} = \sum_{\ell,m=0}^q\delta_{q,\ell+m}\sum_{r=1}^{\ell-1}\sum_{s=1}^{r-1}\delta_{\ell,r+s}(r-s)^2\theta_{e_s,1}\theta_{e_r,1}\Psi^{(m;\{e_1,\ldots,e_s-1,\ldots,e_r-1,\ldots,e_\ell,\ldots,e_q\})}\,,
\end{equation}
and
\begin{equation}
\Omega_1^{(q;\ove)} = \sum_{\ell,m=0}^q\delta_{q,\ell+m}\ell^2\theta_{e_\ell,1}\Psi^{(m;\{e_1,\ldots,e_\ell-1,\ldots,e_q\})}\,.
\end{equation}
In order to complete the proof, we need to show that
$\Omega_0^{(q;\ove)}<\Omega_1^{(q;\ove)}$
$\forall\,\ove\in\cS_q^q$. To this aim, we find it convenient to use
a slightly different representation of $\Psi^{(p;\ove)}$, {\it viz.}
\begin{equation}
\Psi^{(p;\ove)} =
\delta_{\cP_q(\ove),p}\sum_{n=1}^p\sum_{\ovc\in\cS_q^n}{\sum_{k=1}^q
  c_k\choose c_1,\ldots,c_q}{\sum_{k=1}^q (e_k-c_k)\choose
  e_1-c_1,\ldots e_q-c_q}\theta_{\ove,\ovc}\,.
\end{equation}
Suppose that $c_\ell\ge 1$ for some $1\le \ell\le q$. Then,
\begin{align}
\theta_{e_\ell,1}\Psi^{(q-\ell;\{e_1,\ldots,e_\ell-1,\ldots,e_q\})} & \nonumber\\[0.0ex]
& \hskip -3.5cm = \delta_{\cP_q(\ove),q}\sum_{t=1}^q\sum_{\ovc\in\cS_q^t}{\sum_{k=1}^q
  c_k\choose c_1,\ldots,c_q}{\sum_{k=1}^q (e_k-c_k)-1\choose
  e_1-c_1,\ldots,e_n-c_n-1,\ldots
  e_q-c_q}\theta_{e_n,c_n+1}\prod_{s\ne n}^{1\ldots q}\theta_{e_s,c_s}\,.
\end{align}
On performing the change of variable $c_\ell\to c_\ell+1$, we easily arrive
at
\begin{align}
\theta_{e_\ell,1}\Psi^{(q-\ell;\{e_1,\ldots,e_\ell-1,\ldots,e_q\})} =
\delta_{\cP_q(\ove),q}\sum_{t=\ell}^q\sum_{\ovc\in\cS_q^t}
\frac{c_\ell}{\sum_{k=1}^q c_k}{\sum_{k=1}^q
  c_k\choose c_1,\ldots,c_q}{\sum_{k=1}^q (e_k-c_k)\choose
  e_1-c_1,\ldots, e_q-c_q}\theta_{\ove,\ovc}\,.
\label{eq:psin}
\end{align}
Analogously, for $e_s\ge 1$, $e_r\ge1$ and $1\le s<r<q$, we have
\begin{align}
\theta_{e_s,1}\theta_{e_r,1}\Psi^{(q-s-r;\{e_1,\ldots,e_s-1,\ldots,e_r-1,\ldots,e_q\})} & \nonumber\\[0.0ex]
& \hskip -6.2cm =
\delta_{\cP_q(\ove),q}\sum_{t=r+s}^q\sum_{\ovc\in\cS_q^t}
\frac{c_s\,c_r}{\left(\sum_{k=1}^q c_k\right)\left(\sum_{k=1}^q c_k-1\right)}{\sum_{k=1}^q
  c_k\choose c_1,\ldots,c_q}{\sum_{k=1}^q (e_k-c_k)\choose
  e_1-c_1,\ldots,
  e_q-c_q}\theta_{\ove,\ovc}\theta_{\sum_{k=1}^qc_k,2}\,.
\label{eq:psisr}
\end{align}
Note that the lower limit of the sum over $t$ in eq.~(\ref{eq:psin})
can be reduced from $\ell$ to $1$, since $c_\ell=0$ if $\ove\in\cS_q^t$ and
$t<\ell$. The same cannot be done in eq.~(\ref{eq:psisr}) without
increasing the resulting sum. Therefore, 
\begin{equation}
\Omega_1^{(q;\ove)} = \delta_{\cP_q(\ove),q}\sum_{t=1}^q\sum_{\ovc\in\cS_q^t}
\frac{\sum_{\ell=1}^q \ell^2 c_\ell}{\sum_{k=1}^q c_k}{\sum_{k=1}^q
  c_k\choose c_1,\ldots,c_q}{\sum_{k=1}^q (e_k-c_k)\choose
  e_1-c_1,\ldots, e_q-c_q}\theta_{\ove,\ovc}\,,
\end{equation}
and
\begin{align}
\Omega_0^{(q;\ove)} & \le \delta_{\cP_q(\ove),q}\sum_{t=1}^q\sum_{\ovc\in\cS_q^t}
\frac{\sum_{\ell=1}^q\sum_{r=1}^{\ell-1}\sum_{s=1}^{r-1}\delta_{\ell,r+s}(r-s)^2
  c_s\,c_r}{\left(\sum_{k=1}^qc_k\right)\left(\sum_{k=1}^q
  c_k-1\right)}\nonumber\\[2.0ex]
& \hskip 2.6cm \cdot {
  c_k\choose c_1,\ldots,c_q}{\sum_{k=1}^q (e_k-c_k)\choose
  e_1-c_1,\ldots, e_q-c_q}\theta_{\ove,\ovc}\,.
\end{align}
Now, it is immediate to prove that
\begin{equation}
\sum_{\ell=1}^q\sum_{r=1}^{\ell-1}\sum_{s=1}^{r-1}\delta_{\ell,r+s}(r-s)^2
c_s\,c_r<\sum_{\ell=1}^q \ell^2 c_\ell\left(\sum_{k=1}^q
c_k-1\right)\,,\qquad \forall\,\ovc\in\cS_q^t\,.
\end{equation}
Indeed, given $\ovc\in\cS_q^t$ each non--vanishing contribution
$c_sc_r>0$ with $s<r$ is weighted by
$(r-s)^2$ on the \lhs and by $(r^2+s^2)$ on the \rhs The remaining
terms on the \rhs are $\sum_{\ell=1}^{q-1}\ell^2 c_{\ell}(c_\ell-1)\ge
0$ and $q^2c_q(\sum_{k=1}^qc_k-1)\ge 0$. This concludes the proof.
\end{proof}

\section{Concluding remarks}

Conditioning a vector $X\sim\cN_v(0,\Lambda)$ with
$\Lambda=\diag(\lambda)$ to a centered Euclidean ball $\cB_v(\rho)$ of
square radius~$\rho$ affects non--trivially the covariance matrix of
its square components. Since the conditional moments of $X$ cannot be 
calculated in closed--form, the only viable approach (besides
numerical computation) to characterizing the truncational effects
consists in establishing analytic bounds to the conditional
correlations (variances and covariances) of the square
components of $X$. Such estimates are also referred to in the 
literature as {\it square correlation inequalities}. 

In this paper, we specifically focused on the conditional
variances. In particular, our aim was proving
eq.~(\ref{eq:varineq}). The analyses presented in the previous
sections go in this direction, yet they do not solve the problem in a
conclusive way. The arguments proposed apply in the opposite regimes
of strong and weak truncations. For $0<\rho<2\lambda_n$,
eq.~(\ref{eq:varineq}) is easily proved. A bigger effort is required
for $\rho\gg\lambda_n$. Nothing is said regarding the intermediate
region. 
We conclude with two major criticisms, representing at the
same time an outlook of future research:
\begin{itemize}
\item{the weak truncation region is not sharply defined:
  the asymptotic property stated by Theorem~5.1 is
  certainly sufficient to prove that the $p^{\rm th}$ order of the
  weak truncation expansion of $\Delta_n$ is negative at $\rho>\rho_p^*$
  for some $\rho_p^*$, but the theorem does not provide any estimate
  of $\rho_p^*$. A better characterization of the coefficient
  functions $\eta_k$ and $\Xi_{\Delta_n}$ far from the
  asymptotic regime would help identify precise conditions to extend
  the proof of eq.~(\ref{eq:varineq}) to large yet finite values of
  $\rho$ along the same lines of Theorem~5.1;} 
\item{we also lack a general proof of convergence of the weak
  truncation expansion. The argument presented in sect.~4 suggests
  uniform convergence in $v\le 5$ dimensions, but it is based on
  a numerical estimate of the vanishing rate of the $p^{\rm th}$ term of
  the expansion, which cannot be legitimately extrapolated to
  $p\to\infty$.}
\end{itemize}
The weak truncation expansion of a given observable $f$ (built from
Gaussian integrals $\alpha_{k\ell m\ldots}$) is to all extents a perturbative
expansion around the factorized value $f$ takes as
$\rho\to\infty$. As such, it is affected by the usual problems
encountered with perturbative expansions. Having proved a property of
$\Delta_n$ to all orders represents the main (non--trivial)
contribution of the present paper.

\bibliographystyle{plain}
\bibliography{main}

\begin{thebibliography}{1}

\bibitem{abramowitz}
M.~Abramowitz and I.~A. Stegun.
\newblock {\em Handbook of Mathematical Functions with Formulas, Graphs, and
  Mathematical Tables}.
\newblock Dover Publications, New York, 1964.

\bibitem{efron}
B. Efron.
\newblock {Increasing properties of P\'olya frequency functions.}.
\newblock {\em The Annals of Mathematical Statistics}, 36:272-279 1965.

\bibitem{knuth}
R.~L. Graham, D.~E. Knuth, and O.~Patashnik.
\newblock {\em Concrete Mathematics: A Foundation for Computer Science}.
\newblock Addison-Wesley Longman Publishing Co., Inc., Boston, MA, USA, 2nd
  edition, 1994.

\bibitem{joag:1983}
K.~Joag-Dev and F.~Proschan.
\newblock Negative association of random variables with applications.
\newblock {\em The Annals of Statistics}, 11(1):286--295, 1983.

\bibitem{palombi2}
F.~{Palombi} and S.~{Toti}.
\newblock {A perturbative approach to the reconstruction of the eigenvalue
  spectrum of a normal covariance matrix from a spherically truncated
  counterpart}.
\newblock {\em ArXiv e-prints}, July 2012, 1207.1256.

\bibitem{palombi4}
F.~{Palombi}, S.~{Toti}, and R.~{Filippini}.
\newblock {Numerical reconstruction of the covariance matrix of a spherically
  truncated multinormal distribution}.
\newblock {\em ArXiv e-prints}, February 2012, 1202.1838.

\bibitem{palombi1}
F.~Palombi, S.~Toti, R.~Filippini and V. Tomeo.
\newblock {A forward search algorithm for compositional data}
\newblock {\em Key Invited Paper, Conference of European Statisticians, 
  Work Session on Statistical Data Editing, Ljubljana, Slovenia, 9--11 May 2011.}
\newblock {\scriptsize http://www.unece.org/fileadmin/DAM/stats/documents/ece/ces/ge.44/2011/wp.41.e.pdf}

\bibitem{prekopa}
A.~Pr\'ekopa.
\newblock Logarithmic concave measures with applications to stochastic
  programming.
\newblock {\em Acta Scientiarum Mathematicarum}, 32:301--316, June 1971.

\bibitem{ruben}
H.~Ruben.
\newblock Probability content of regions under spherical normal distributions,
  I.
\newblock {\em The Annals of Mathematical Statistics}, 31(3):598--618,
  September 1960.




\end{thebibliography}

\end{document}